\documentclass[12pt,leqno]{article}
\usepackage{amsmath,amsthm,amssymb}
\usepackage[curve,matrix,arrow]{xy}
\newcommand{\smallbullet}{{\mathchoice
{\raise.8pt\hbox{$\scriptstyle\bullet$}}
{\raise.8pt\hbox{$\scriptstyle\bullet$}}
{\raise.6pt\hbox{$\scriptscriptstyle\bullet$}}
{\hbox{$\scriptscriptstyle\bullet$}}}}
\newcounter{x}\newenvironment{romanlist}{
 \begin{list}{\rm\hskip-2.0em\hbox to 2.0em{\hfil(\roman{x})}}
 {\usecounter{x}
 \setlength{\leftmargin}{2.5em}\setlength{\itemsep}{6pt}
 \setlength{\parsep}{0pt}}}
 {\end{list}}
\newcounter{xx}\newenvironment{romanilist}{
 \begin{list}{\rm\ifnum\arabic{xx}=1\hskip-2.0em\hbox to
 2.0em{\hfil(\roman{xx})\hskip.65em}\else\hskip-2.0em\hbox to
 2.0em{\hfil(\roman{xx})}\fi}
 {\usecounter{xx}
 \setlength{\leftmargin}{2.5em}\setlength{\itemsep}{6pt}
 \setlength{\parsep}{0pt}}}
 {\end{list}}
\newcounter{xxx}\newenvironment{romanslist}{
 \begin{list}{\rm\ifnum\arabic{xxx}=1\hskip-2.0em\hbox to
 11.7em{\hfil(\cite[Proposition 3.1]{f})\hspace{.7em}(\roman{xxx})\hskip.65em}\else\hskip-2.0em\hbox to
 2.0em{\hfil(\roman{xxx})}\fi}
 {\usecounter{xxx}
 \setlength{\leftmargin}{2.5em}\setlength{\itemsep}{6pt}
 \setlength{\parsep}{0pt}}}
 {\end{list}}
\newcommand{\firstitem}{\item \hskip-.65em}
\newtheorem{thm}{Theorem}[subsection]
\newtheorem{pr}[thm]{Proposition}
\newtheorem{lm}[thm]{Lemma}
\newtheorem{cor}[thm]{Corollary}
\theoremstyle{definition}
\newtheorem{defn}[thm]{Definition}
\theoremstyle{remark}
\newtheorem{rem}[thm]{Remark}
\newtheorem{conv}[thm]{Convention}
\newtheorem{nota}[thm]{Notation}
\numberwithin{equation}{subsection}
\DeclareMathOperator{\Hom}{Hom}
\DeclareMathOperator{\Bl}{B\ell}
\DeclareMathOperator{\Spec}{Spec}

\DeclareMathOperator{\Sym}{Sym}

\DeclareMathOperator{\intrior}{int}

\DeclareMathOperator{\Isom}{Isom}
\DeclareMathOperator{\VVert}{Vert}
\DeclareMathOperator{\rk}{rk}
\DeclareMathOperator{\im}{im}
\newcommand{\bP}{{\mathbb P}}
\newcommand{\A}{{\mathbb A}}

\newcommand{\Q}{{\mathbb Q}}
\newcommand{\R}{{\mathbb R}}
\newcommand{\Z}{{\mathbb Z}}
\newcommand{\cO}{{\mathcal O}}
\newcommand{\cE}{{\mathcal E}}
\newcommand{\cF}{{\mathcal F}}
\newcommand{\mm}{{\mathfrak m}}
\newcommand{\Aff}{{\it Aff}}
\newcommand{\smallcap}{\mathbin{\raise.4pt\hbox{$\scriptstyle\cap$}}}
\newcommand{\smallcirc}{\mathchoice
{\mathbin{\raise.9pt\hbox{$\scriptstyle\circ$}}}
{\mathbin{\raise.9pt\hbox{$\scriptstyle\circ$}}}
{\mathbin{\raise.35pt\hbox{$\scriptscriptstyle\circ$}}}
{\mathbin{\hbox{$\scriptscriptstyle\circ$}}}}
\newcommand{\smallamalg}{\mathchoice
{\mathbin{\raise.6pt\hbox{$\scriptstyle\amalg$}}}
{\mathbin{\raise.6pt\hbox{$\scriptstyle\amalg$}}}
{\mathbin{\raise.35pt\hbox{$\scriptscriptstyle\amalg$}}}
{\mathbin{\hbox{$\scriptscriptstyle\amalg$}}}}
\newcommand\eqto{\stackrel{\lower1.5pt\hbox{$\scriptstyle\sim\,$}}\to}

\begin{document}
\title{Cycle groups for Artin stacks}
\author{Andrew Kresch$^1$}
\date{28 October 1998}
\maketitle
\footnotetext[1]{Funded by a
Fannie and John Hertz Foundation Fellowship
for Graduate Study and an Alfred P. Sloan Foundation Dissertation
Fellowship}

\tableofcontents

\section{Introduction}

We define a Chow homology functor $A_*$ for Artin stacks and prove that it
satisfies some of the basic properties expected from intersection theory.
Consequences include an integer-valued intersection product
on smooth Deligne-Mumford stacks, an affirmative answer
to the conjecture that any smooth stack with finite but possibly nonreduced
point stabilizers should
possess an intersection product (this provides a positive answer
to Conjecture 6.6 of \cite{v}), and more generally an intersection product
(also integer-valued) on smooth Artin stacks
which admit stratifications by global quotient stacks.

The definition presented here generalizes existing definitions.
For representable stacks (i.e., algebraic spaces), the functor $A_*$
defined here
reproduces the standard Chow groups.
For Deligne-Mumford stacks \cite{dm} the functor differs by torsion
from the na\"\i{}ve Chow groups
(algebraic cycles modulo rational equivalence).
The na\"\i{}ve Chow groups lead to a $\Q$-valued intersection
theory on Deligne-Mumford stacks \cite{g,v}.
However, there is evidence that the na\"\i{}ve Chow group functor
is not the correct object to work with if we wish to have
integer coefficients.
For instance, with the na\"\i{}ve Chow groups,
there is no integer-valued intersection product
on smooth Deligne-Mumford stacks.
Even Chern classes of vector bundles do not exist except with
rational coefficients.

Many Deligne-Mumford stacks are also global quotient stacks,
and for a general global quotient stack,
the functor $A_*$ reproduces the
equivariant Chow groups of Dan Edidin and William Graham \cite{eg}.
In fact, the realization that the equivariant Chow groups are the
correct groups to work with
provided the starting point for the definition presented here.
The main idea of \cite{eg} is that a global quotient stack possesses
a vector bundle $E$ such that
the total spaces of $E^{\oplus n}$
(sums of $n$ copies of $E$) become increasingly well approximated
by algebraic spaces as $n$ gets large.
In any fixed codimension,
the (na\"\i{}ve) Chow groups of these approximating spaces stabilize, and
serve as the equivariant Chow homology functor.

This preprint takes the idea further by considering
cycle classes in {\em all} vector bundles on an Artin stack.
We say that a bundle $F$ dominates
$E$ if there exists a vector bundle surjection $\varphi\colon F\to E$.
Once we know that corresponding pullback map on cycle classes is
independent of the particular choice of $\varphi$,
we can take the direct limit of over all vector bundles.
This gives us, for global quotient stacks, the equivariant Chow groups.
In general, we obtain groups $\widehat A_*X$, but
these groups themselves do not provide intersection theory.
One can produce a stack $X$ which has no nontrivial
vector bundles, but which has a closed substack $Y$
which is a global quotient.
There is no way to push forward a nontrivial class
$\alpha\in\widehat A_*Y$
via the inclusion map $i\colon Y\to X$.
The only thing to do is to
define $i_*\alpha$ to be the formal pair $(i,\alpha)$.
The set of all such pairs, modulo what is more or less the
weakest possible equivalence relation which guarantees functoriality of
pushforward, forms a group $A_*X$ (section \ref{homologyfunctor}) which
provides a reasonable amount of intersection theory.

Once we have proved some basic properties for $A_*$
(sections \ref{basic}--\ref{excsection} and section \ref{vb}),
we can obtain
an integer-valued intersection product on
smooth Deligne-Mumford stacks (section \ref{elementary}).
All we need for this are standard constructions
and properties from intersection theory.
By way of constrast, intersection theory on Artin stacks
requires homotopy invariance for objects more general
than vector bundles.
These objects, called {\em vector bundle stacks} in \cite{bf},
look locally like a quotient of of one vector bundle
by the (additive) action of another.
The proof of homotopy invariance for vector bundle stacks
uses a localization argument
(section \ref{xxa}),
and once we have this property for a class of stacks
then intersection theory follows.
The relevant class of stacks consists of stacks which can be
stratified by locally closed substacks that are global quotients.
So, any such stack, if it is smooth, possesses an integer-valued
intersection product (section \ref{intartsta}).
This class of stacks includes (finite-type approximations of) many interesting
stacks such as moduli stacks of stable or pre-stable curves
and moduli stacks of vector bundles.

The paper \cite{bf} constructs a virtual fundamental class of the
expected dimension
from a perfect obstruction theory.
This class is gotten as the
``intersection with the zero section'' of a cone stack sitting
in some vector bundle stack.
Lacking intersection theory on Artin stacks
(a vector bundle stack is an Artin stack),
the authors were forced to impose a technical hypothesis in
order to carry out their construction.
The construction can now be done in general
(section \ref{virfundclass}).
We conclude this preprint with a discussion (section \ref{localiz}) of the
localization formula for torus actions.

The author is greatly indebted to his thesis advisor William Fulton,
who provided constant guidance and assistance during the
course of this project.
The author would like to acknowledge helpful discussions with
Dan Edidin, Henri Gillet, and Charles Walter.
Much of the work was done at the Mittag-Leffler Institute during the
1996--97 program in algebraic geometry.
The author would like to thank the staff and organizers there for
providing a fantastic research environment.

\section{Definition and first properties}
\label{firstproperties}

\subsection{The homology functor}
\label{homologyfunctor}

We start by recalling the notion of algebraic cycle and
rational equivalence on a stack and by introducing the functor $A^\circ_*$
of cycles modulo rational equivalence.  From
this we build the functor $A_*$ by a succession of direct limits.
The first limit is over vector bundles (Definition \ref{setby}, below).
The second limit is over projective morphisms to the target stack
(Definition \ref{setay} (ii)).
As is always the case when we wish to take a limit of abelian groups
indexed by a directed set which comes from a
category by collapsing morphisms, we must take care to verify
that the induced map on groups is independent of the choice of morphism
whenever more than one morphism exists between two objects of the category
(Remark \ref{twovbsurj} for the limit over vector bundles
and Remark \ref{twoyincl} for the limit over projective morphisms).

\begin{conv}
\label{convone}
All stacks are Artin stacks (i.e., algebraic stacks
with smooth atlases \cite{a,l}) and are
of finite type over a fixed base field.
All morphisms are morphisms over the base field.
All regular (local) immersions are of constant codimension
(a local immersion is by definition a representable unramified morphism;
to every local immersion there is an associated normal cone, and if
the cone is a bundle then the morphism is called a
regular local immersion \cite{v}).
A morphism of stacks $X\to Y$ is called projective if it can
be factored (up to 2-isomorphism)
as a closed immersion followed by the projection
morphism ${\mathbf P}(\cE)\to Y$ coming from
a coherent sheaf $\mathcal E$ of $\cO_Y$-modules on $Y$
\cite{ega, l}.
Closed immersions and projections of the form
${\mathbf P}(\cE)\to Y$ are examples of representable morphisms, so that
every projective morphism of stacks is representable as well.
\end{conv}

Let us recall that if $Y$ is a stack, then there is a category,
known as the {\em category of $Y$-stacks}, whose objects consist
of pairs $(X,f)$ with $X$ a stack (in our case, according to
Convention \ref{convone}, always of finite type over the base field) and
$f$ a representable morphism from $X$ to $Y$.
A morphism from $(X,f)$ to $(X',f')$ is a pair $(\varphi,\alpha)$,
where $\varphi$ is a morphism from $X$ to $X'$, and
$\alpha$ is a 2-morphism from $f$ to $f'\smallcirc\varphi$.
The full subcategory consisting of all $(X,f)$ such that $f$ is projective
is called the {\em category of projective $Y$-stacks}.

\begin{defn}
\label{setay}
\begin{romanilist}
\firstitem An {\em inclusion of components}
is a morphism
$f\colon X\to X'$ which is an isomorphism of $X$ onto
a union of connected components of $X'$.
If $Y$ is a stack, $(X,f)$ and $(X',f')$ are $Y$-stacks, and $g$
is a morphism from $(X,f)$ to $(X',f')$ in the category of
$Y$-stacks (in short, a $Y$-morphism),
then we say that $g$ is a {\em $Y$-inclusion of components}.
The set of isomorphism classes of $Y$-stacks forms a directed set
with $(X,f)\preceq(X',f')$ whenever there exists
a morphism from $(X,f)$ to $(X',f')$
which is a $Y$-inclusion of components.
\item Let $Y$ be a stack.
We denote by ${\mathfrak A}_Y$ the set of isomorphism classes
of projective $Y$-stacks with the partial ordering of (i).
For each $Y$, ${\mathfrak A}_Y$ is a directed set.
\end{romanilist}
\end{defn}

\begin{defn}
\label{setby}
Let $Y$ be a stack.
We denote by ${\mathfrak B}_Y$ the set of isomorphism classes of
vector bundles over $Y$,
partially ordered by declaring $E\preceq F$ whenever there exists a
surjection of vector bundles $F\to E$.
\end{defn}

\begin{defn}
\label{bigdefn}
\begin{romanilist}
\firstitem For $Y$ a stack,
$Z_*Y$ denotes the group of
{\em algebraic cycles} on $Y$, i.e., the free abelian
group on the set of integral closed substacks of $Y$, graded by dimension
\cite{g,v}.
We denote by $W_*Y$ the group of {\em rational equivalences}
on $Y$.
If $k(Z)^*$ denotes the multiplicative
group of rational functions on an integral
substack $Z$, not identically zero,
then $W_jY$ is the direct sum of $k(Z)^*$ over integral closed
substacks $Z$ of dimension $j+1$.
There is a map $\partial\colon W_jY\to Z_jY$ which locally for
the smooth topology sends a rational function to the corresponding
Weil divisor.
\item
The {\em na\"\i{}ve Chow groups} of $Y$
are defined to be $A^\circ_kY=Z_kY/\partial W_kY$.
\item
The {\em Edidin-Graham-Totaro Chow groups}$^1$\footnotetext[1]{In \cite{eg}
the authors attribute the idea behind their construction to
Burt Totaro.} \cite{eg} are defined, for $Y$
connected, by
$\widehat A_kY=\varinjlim_{{\mathfrak B}_Y}
A^\circ_{k+\rk E} E$, and for
$Y=Y_1\amalg \cdots \amalg Y_r$ with each $Y_i$ connected, by
$\widehat A_kY=\bigoplus_{i=1}^r \widehat A_k Y_i$.
\item
Given a morphism $f\colon X\to Y$ with $X$ connected, we define
the {\em restricted Edidin-Graham-Totaro Chow groups} to be the groups
$\widehat A^f_k X = \varinjlim_{{\mathfrak B}_Y}
A^\circ_{k+\rk E}f^*E$.
If
$X=X_1\amalg \cdots \amalg X_r$ with each $X_i$ connected, we set
$\widehat A^f_kX=\bigoplus_{i=1}^r \widehat A^f_k X_i$.
There is a natural map $\iota_f\colon \widehat A^f_k X\to \widehat A_kX$.
\item
If $f\colon X\to Y$ is a projective morphism,
the {\em restricted projective pushforward} is the
map $f_*\colon \widehat A^f_kX\to \widehat A_kY$ defined by
each $\tilde f_*\colon A^\circ_kf^*E\to A^\circ_kE$
(for every $E$ we are letting
$\tilde f\colon f^*E\to E$ denote the pullback of $f$).
\end{romanilist}
\end{defn}

\begin{rem}
\label{twovbsurj}
If $E$ and $F$ are vector bundles on a stack $Y$,
then any two vector bundle surjections $\varphi,\psi\colon E\to F$
induce the same map on na\"\i{}ve Chow groups.
Indeed, let $\Psi_t=\varphi+t(\psi-\varphi)$ for $t$ in the base field.
This defines a morphism $\Psi\colon E\times\A^1\to F$.
Given any closed integral substack $Z$ of $F$, the closure
$\overline{\Psi^{-1}(Z)}$ of $\Psi^{-1}(Z)$
in $E\times\bP^1$ exhibits a rational equivalence between
$\varphi^*[Z]$ and $\psi^*[Z]$ (the rational equivalence between
the fiber over $t=0$ and the fiber over $t=1$ pushes forward to $E$).
Thus $\{A^\circ_{r+\rk E}E\}$ forms a direct system of abelian groups
over ${\mathfrak B}_Y$.
\end{rem}

\begin{rem}
\label{ahatspace}
The natural map
$A^\circ_*Y\to \widehat A_*Y$ is an isomorphism for any
scheme $Y$ \cite[Theorem 3.3, (a)]{f},
or more generally for any algebraic space \cite{k} $Y$
(any algebraic space has a dense open subspace represented by
a scheme, so projection from a vector bundle induces
a pullback map on $A^\circ_*$ which is surjective by the
argument of \cite[Proposition 1.9]{f},
and injectivity is demonstrated exactly as for schemes).
\end{rem}

\begin{rem}
The groups $\widehat A_*Y$ are defined in \cite{eg} only
in the special case of a quotient stack
$Y\simeq[X/G]$ (with $G$ an algebraic group acting on an
algebraic space $X$; the authors employ the notation
$A^G_*X$ but point out that these groups are in fact
invariants of the underlying stack $Y$).
Because quotient stacks admit suitable approximations by
algebraic spaces, the limit in Definition \ref{bigdefn} (iii)
stabilizes after some point,
and now because the pullback map on Chow groups induced
by a vector bundle over an algebraic space is an isomorphism,
the content of Remark \ref{twovbsurj} becomes trivial and therefore
not an ingredient in the construction of \cite{eg}.
\end{rem}

\begin{rem}
An inclusion of components gives rise to unrestricted
projective pushforward, i.e., if $f\colon X\to Y$ is an inclusion
of components then $\iota_f$ is an isomorphism.
\end{rem}

\begin{defn}
\label{bhat}
Let $X$ be a stack.  If $T$ is a stack, and $p_1$ and $p_2$ are
projective morphisms $T\to X$, then
the set
$$\{\,p_{2{*}}\beta_2-p_{1{*}}\beta_1   \mid
(\beta_1,\beta_2)\in \widehat A^{p_1}_kT\oplus \widehat A^{p_2}_kT
{\rm\ satisfies\ }
\iota_{p_1}(\beta_1)=\iota_{p_2}(\beta_2)\,\}$$
is a subgroup of $\widehat A_kX$ which we denote
$\widehat B_k^{p_1,p_2}X$.
If $(X,f)$ is a $Y$-stack,
the union of the subgroups $\widehat B_k^{p_1,p_2}X$, over all $T$ and all
pairs of projective morphisms $p_1$ and $p_2$ such that
$f\smallcirc p_1$ is 2-isomorphic to $f\smallcirc p_2$,
is a subgroup of $\widehat A_kX$ which
we denote $\widehat B_kX$.
\end{defn}

\begin{rem}
\label{twoyincl}
Let $Y$ be a stack, and suppose we are given two $Y$-inclusions of components
$f_1,f_2\colon X\to X'$.
Then for all $\alpha\in\widehat A_*X$,
we have $f_{2{*}}\alpha-f_{1{*}}\alpha\in\widehat B^{f_1,f_2}_k X$
(take $\beta_1=\beta_2=\alpha$).
\end{rem}

\begin{defn}
\label{defnapaj}
Let $Y$ be a stack.
We define
$$A_kY = \varinjlim_{{\mathfrak A}_Y}
( \widehat A_kX / \widehat B_kX ).$$
\end{defn}

The main result of this preprint is

\begin{thm}
\label{mainthm}
Let $k$ be a base field.
The functor $A_*$ of Definition \ref{defnapaj},
from the category of Artin stacks of
finite type over $k$ to the category of
abelian groups graded by dimension, is
contravariant for morphisms which are flat of locally constant relative
dimension, covariant for projective morphisms, and is related
to the functors $A^\circ_*$ and $\widehat A_*$ via
natural maps
$$A^\circ_*X \to \widehat A_* X \to A_*X.$$
There is a canonically defined ring structure
on $A_*X$ when $X$ is smooth and
can be stratified by locally closed substacks which are each
isomorphic to the quotient stack of an algebraic group acting on an
algebraic space.
The functor $A_*$ satisfies the following properties:
\begin{romanlist}
\item For any algebraic space $X$, the map
$A^\circ _*X\to A_*X$ is an isomorphism of groups,
and for $X$ smooth is an isomorphism of rings.
\item For any Deligne-Mumford stack $X$,
the map $A^\circ _*X\otimes\Q\to A_*X\otimes\Q$ is an isomorphism of groups,
and for $X$ smooth is an isomorphism of rings.
\item For any algebraic space $X$ with action of a
linear algebraic group $G$,
the map $\widehat A_*[X/G] \to A_*[X/G]$
is an isomorphism of groups,
and for $X$ smooth is an isomorphism of rings.
\item For any stack $X$ with closed substack $Z$ and complement $U$,
the excision sequence
$A_j Z\to A_j X\to A_j U\to 0$
is exact.
\item We have $A_jX=0$ for all $j>\dim X$.
\item If $\pi\colon E\to X$ is a vector bundle of rank $e$,
then the induced pullback map $\pi^*\colon A_j X\to A_{j+e} E$ is
an isomorphism.
\item If $\pi\colon E\to X$ is a vector bundle of rank $e$
with associated projective bundle $p\colon P(E)\to X$ and
line bundle $\cO_E(1)$ on $P(E)$, then the map
$$\theta_E\colon \bigoplus_{i=0}^{e-1} A_{j-e+i+1} X\to A_j P(E)$$
given by $(\alpha_i)\mapsto \sum_{i=0}^{e-1}
c_1(\cO_E(1))^i\smallcap p^*\alpha_i$
is an isomorphism.
\item There are Segre classes and Chern classes of vector bundles, and these
satisfy the usual universal identities.
\item There are Gysin maps for regular immersions and
regular local immersions, and these maps are
functorial, commute with each other, and are compatible with
flat pullback and projective pushforward.
\item If $X$ can be stratified by locally closed substacks which are
isomorphic to quotient stacks,
then for any vector bundle stack $\pi\colon B\to X$
of virtual rank $e$,
the induced pullback map $\pi^*\colon A_jX\to A_{j+e}B$ is
an isomorphism.
\item Gysin maps exist for l.c.i.\ morphisms $f\colon X\to Y$
whenever $X$ can be stratified by quotient stacks,
and these maps are
functorial, commute with each other, and are compatible with
flat pullback and projective pushforward.
\end{romanlist}
\end{thm}

\begin{rem}
As we shall see,
the class of stacks which can be stratified by quotient stacks
includes Deligne-Mumford stacks, and more generally
stacks with quasifinite diagonal (Proposition \ref{isaquotient}).
This class of stacks is stable under representable morphisms
and formation of products (Proposition \ref{gq}).
\end{rem}

\begin{nota}
We denote a typical element of $A_kY$ by
$(f,\alpha)$ with $f\colon X\to Y$ projective
and $\alpha\in\widehat A_kX$.
When we speak of a {\em cycle} we refer to such a choice
among all the representatives of a given cycle class.
Unless specifically stated to the contrary,
an {\em identity of cycles} refers to
an identity in $\widehat A_*$ of the relevant stack.
\end{nota}

\begin{rem}
\label{equivzero}
Explicitly, now, if $(X,f)$ is a projective $Y$-stack,
an element $(f,\alpha)$ of $A_kY$
is equivalent to zero if and only if there exists a
$Y$-inclusion of components $i\colon X\to X'$ for some
projective $Y$-stack $(X',f')$
such that there exist projective morphisms
$p_1,p_2\colon T\to X'$ and
$\beta_i\in A^{p_i}_kT$ ($i=1,2$) such that
$f'\smallcirc p_1$ is 2-isomorphic to $f'\smallcirc p_2$, we have
$\iota_{p_1}(\beta_1)=\iota_{p_2}(\beta_2)$ in $\widehat A_*T$,
and $i_*\alpha=p_{2{*}}\beta_2-p_{1{*}}\beta_1$ in $\widehat A_*X'$.
\end{rem}

\begin{rem}
\label{almostequiv}
Suppose $f\colon X\to Y$ and $p_1,p_2\colon T\to X$ are
projective and $g:=f\smallcirc p_1$ and $f\smallcirc p_2$ are 2-isomorphic.
Then, for any $\beta_1\in\widehat A^{p_1}_*T$ and
$\beta_2\in\widehat A^{p_2}_*T$,
we have $(f,p_{2{*}}\beta_2-p_{1{*}}\beta_1)=
(g,\iota_{p_2}(\beta_2)-\iota_{p_1}(\beta_1))$
in $A_*Y$
(consider $q_1,q_2\colon T\amalg T\to X\amalg T$ given by
$q_1=p_1\smallamalg 1_T$ and $q_2=1_T\smallamalg p_2$).
In particular, $(g,\iota_{p_1}(\beta_1))=(f,p_{1{*}}\beta)$
in $A_*Y$.
This plus Remarks \ref{ahatspace} and
\ref{equivzero} establishes that for
any algebraic space $Y$, the natural map
$A^\circ_*Y\to A_*Y$ is an isomorphism.
Similarly, from the fact that whenever $\pi\colon E\to X$ is
a vector bundle over a Deligne-Mumford stack $X$, the
pullback map $\pi^*$ induces an isomorphism
$A^\circ_jX\otimes\Q\to A^\circ_{j+\rk E}E\otimes\Q$,
we conclude that the natural map
$A^\circ_*Y\otimes\Q\to A_*Y\otimes\Q$ is an isomorphism.
\end{rem}

\begin{rem}
\label{linalggpacts}
Let a linear algebraic group $G$ act on an algebraic space $V$,
and let $Y=[V/G]$ be the stack quotient of $V$ under the action of $G$.
Then $A_jY$ is the equivariant Chow group $A^G_jV$ of \cite{eg}.
Indeed,
suppose $f\colon X\to Y$ is projective and $E\to X$ is a vector bundle.
We consider a $r$-dimensional
representation of $G$ such that the
corresponding action of $G$ on affine space is free off of some
locus of codimension larger than $\dim X - j$
(such exists for suitable $r$).
Such a representation corresponds to a vector bundle on $BG$;
we let $F$ denote the pullback of this bundle to $Y$ and
$F'$ the pullback to $X$.
Now the map $A^\circ_{j+r}F'\to A^\circ_{j+r+\rk E} E\oplus F'$ induced by
pullback is an isomorphism,
so the map $A^\circ_{j+\rk E} E\to \widehat A_j X$
factors through $A^\circ_{j+r}F'$,
and now for any $\alpha\in A^\circ_{j+\rk E} E$ if we let
$\alpha'$ denote the element of $A^\circ_{j+r}F'$
determined by this factorization,
then we have $(f,\alpha)=(1_Y,f'_*\alpha')$,
where $f'\colon F'\to F$ denote the pullback of $f$.
The resulting maps $A^\circ_{j+\rk E}E\to A^G_jV$ determine a map
$A_jY\to A^G_jV$ which is then seen to be inverse to the
natural map $A^G_jV \simeq \widehat A_jY\to A_jY$.
\end{rem}

\subsection{Basic operations on Chow groups}
\label{basic}

The basic operations we consider here are flat pullback,
projective pushforward, and Gysin maps for principal effective
Cartier divisors.

Given a flat morphism $f\colon Y'\to Y$ of locally constant relative
dimension (but not necessarily representable),
we define a pullback operation $f^*$ as follows.
First, we observe that if $h\colon X\to Y$ is projective
and if we form the fiber diagram
$$\xymatrix
{ X' \ar[r]^{f'} \ar[d]_{h'} & X \ar[d]^h \\
Y' \ar[r]^f & Y
}$$
then $f'$ is also flat of locally constant relative dimension,
and if for any vector bundle $E\to X$ we denote by $\tilde f'$ the
pullback of $f'$ via $E\to X$ then the maps
$\tilde f^{\prime{*}}\colon A^\circ_*E\to A^\circ_*E'$
define a map
$f^{\prime{*}}\colon \widehat A_*X\to \widehat A_*X'$
(there is, of course, a shift in grading by the local relative
dimensions of $f$).
Since flat pullback on $A^\circ_*$ commutes with projective pushforward,
the map $f^{\prime{*}}$ descends to give a map
$\widehat A_*X/\widehat B_*X\to\widehat A_*X'/\widehat B_*X'$
and passes to the limit to give us
$$f^*\colon A_*Y\to A_*Y'.$$

As hinted in the introduction, the definition of proper pushforward
is a tautology.
Given a projective morphism $g\colon Y\to Z$, we define
$g_*(f,\alpha)=(g\smallcirc f,\alpha)$.
By the remarks above,
this map is compatible with restricted projective pushforward:
$g_*(1_Y,\iota_g(\alpha))=(1_Y,g_*\alpha)$.
More generally, suppose $h\colon W\to Z$ is projective.
If we form the fiber square
$$\xymatrix{
X \ar[r]^{g'} \ar[d]_{h'} & W \ar[d]^h \\
Y \ar[r]^g & Z}$$
and if $\alpha\in \widehat A^{g'}_kX$, then we have
$g_*(h',\iota_{g'}(\alpha))=(h,g'_*\alpha)$.
In a fiber diagram, projective pushforward on $A_*$ commutes with
flat pullback.
This follows in a routine fashion from the similar fact for $A^\circ_*$.

Suppose $Y$ is a stack and we are given a morphism $\varphi\colon Y\to\A^1$.
We denote by $Y_0$ the fiber of $\varphi$ over 0 and by
$s$ the inclusion $Y_0\to Y$.
There is a morphism $s^*\colon Z_kY\to Z_{k-1}Y_0$ which sends
$[V]$ to $[(\varphi|_V)^{-1}(0)]$ if $V\not\subset Y_0$
and sends $[V]$ to 0 otherwise
\cite[Remark 2.3]{f}, as well as a compatible morphism
$s^*\colon W_kY\to W_{k-1}Y_0$ which sends a rational function $r$
on $V\not\subset Y_0$ to the formal sum of all the components of the
tame symbol of $r$ and $\varphi$ which are supported in $Y_0$
\cite{glms,kr}.
Thus we obtain
a map on rational equivalence classes of cycles
$s^*\colon A^\circ_kY\to A^\circ_{k-1}Y_0$.
This map respects flat pullback, so we obtain
$s^*\colon \widehat A_kY\to \widehat A_{k-1}Y_0.$

By compatibility with proper pushforward we have the following
form of the projection formula:
\begin{pr}
\label{egtprform}
Let $Y$ be a stack with morphism to $\A^1$ and let $f\colon X\to Y$ be
a projective morphism.
If we form the fiber diagram
$$\xymatrix{
X_0 \ar[r]^{f_0} \ar[d]_t & Y_0 \ar[r] \ar[d]_s & {\{0\}} \ar[d] \\
X \ar[r]^f & Y \ar[r] & {\A^1}
}$$
then
$f_{0{*}}(t^*\alpha)=s^*(f_*\alpha)$ for
all $\alpha\in \widehat A^f_kX$
(where $f_*$ and $f_{0{*}}$ denote restricted projective pushforward).
\end{pr}
Given such a fiber diagram,
the map $(f,\alpha)\mapsto (f_0,t^*\alpha)$ specifies the {\em Gysin map}
$s^* \colon A_kY\to A_{k-1}Y_0$.
That this map respects equivalence is an application of
Proposition \ref{egtprform}, and
is left to the reader.
More generally, if $s\colon Y_0\to Y$ is a closed immersion
and if there is a neighborhood $U$ of $Y_0$ and a
function $\varphi\colon U\to \A^1$ such that
$Y_0=U\times_{\A^1}\{0\}$,
then the composite $A_kY\to A_kU\to A_{k-1}Y_0$ is
independent of $U$ and we call this map the Gysin map $s^*$.
For instance, a morphism $Y\to\bP^1$ determines a Gysin map
$A_kY\to A_{k-1}(Y\times_{\bP^1}\{0\})$.

\subsection{Results on sheaves and vector bundles}

This section collects several elementary facts.

\begin{pr}
\label{xtendcoh}
Let $\cE$ be a coherent
sheaf on an open substack $U$ of a stack $X$.
Then there exists a coherent sheaf $\cE'$ on $X$
such that $\cE'|_U\simeq\cE$.
\end{pr}

\begin{proof}
This is \cite[Corollary 8.6]{l}.
\end{proof}

\begin{cor}
\label{xtendproj}
Let $X$ be a stack and let $U$ be an open substack.
Given a projective morphism $g\colon S\to U$ with $S$ reduced, there exists
a fiber diagram
$$\xymatrix{
S\ar[r]\ar[d]^g & T\ar[d]^f \\
U\ar[r] & X
}$$
with $f$ projective.
\end{cor}

\begin{pr}
\label{xtend}
Let $X$ be a stack, with $U$ an open substack.
Given any vector bundle $E\to U$, there exists a projective
morphism $X'\to X$ which is an isomorphism restricted to $U$,
and a vector bundle $E'\to X'$, such that the restriction of
$E'$ to $U$ is isomorphic to $E$.
\end{pr}

\begin{proof}
Let $r$ be the rank of $E$.
Let $\cE$ be the sheaf of sections of $E$.  We may extend $\cE$ to
a coherent sheaf on all of $X$.
If $X'\to X$ is the Grassmannian of rank $r$ locally
free quotients of $\cE$ and $E'$ is the universal quotient bundle,
then $X'\to X$ is an isomorphism over the locus where $\cE$ is
locally free (in particular, over $U$), and $E'$ restricted to this
locus agrees with the vector bundle determined by $\cE$.
There exists a closed immersion $X'\to{\mathbf P}(\bigwedge^r\cE)$ which
exhibits $X'\to X$ as a projective morphism.
\end{proof}

\begin{pr}
\label{xtendxtension}
Let $X$ be a stack, and let $U$ be the complement of a Cartier
divisor $D$ on $X$.
Let $E$ and $F$ be vector bundles on $X$, and suppose
an extension
\begin{equation}
\label{xtension}
0\to F|_U\to Q\to E|_U\to 0
\end{equation}
on $U$ is given.
Then there exists an extension of $E$ by $F\otimes\cO(nD)$ for some $n$
which restricts to the extension (\ref{xtension}) on $U$.
\end{pr}

\begin{proof}
We need only take $n$ sufficiently large so that
the extension class of $Q$ lies in the image of
$H^1(X,E^*\otimes F\otimes \cO(nD))\to H^1(U,(E^*\otimes F)|_U)$.
\end{proof}

In the statement of Proposition \ref{xtendxtension},
the term ``Cartier divisor'' refers to a substack of codimension 1
which on some smooth atlas is given by the vanishing of a single
function which is a non zero divisor.
Since blowing up exhibits any open substack as the
complement of a Cartier divisor, we have

\begin{cor}
\label{xtendxtensioncor}
Let $X$ be a stack, and let $U$ be an open substack.
If vector bundles $E$ on $X$ and $F_0$ on $U$ and an extension
$0\to F_0\to Q_0\to E|_U\to 0$ are given,
then there exist a projective morphism $X'\to X$ which is
an isomorphism over $U$, a vector bundle $F$ on $X$ such that
$F|_U\simeq F_0$, and an extension
$0\to F\to Q\to E\to 0$ on $X'$ whose restriction to $U$ is the
given extension.
\end{cor}

\subsection{The excision sequence}
\label{excsection}

\begin{pr}
\label{exci}
Let $Y$ be a closed substack of $X$ with
inclusion map $\sigma$.
Let $U$ be the complement of $Y$ in $X$, with $\rho\colon U\to X$.
Then the sequence
$$A_*Y\stackrel{\sigma_*}\to A_*X \stackrel{\rho^*}\to A_*U\to 0$$
is exact.
\end{pr}

\begin{proof}
For surjectivity on the right, suppose $(f,\alpha)$ is a
cycle in $A_*U$, with $f\colon S\to U$ projective.
By Corollary \ref{xtendproj}, $f$ is the
restriction to $U$ of some projective morphism $T\to X$.
If $\alpha$ is represented by a cycle in a vector bundle $E\to S$,
then by Proposition \ref{xtendcoh}, there exists $T'\to T$
projective, which is an isomorphism over $S$, and a bundle
$E'\to T$ such that $E'|_S\simeq E$.
Now the desired result follows since the restriction map
$A^\circ_*E'\to A^\circ_*E$ is surjective.

Since the composite $\rho^*\smallcirc\sigma_*$ is clearly zero,
all that remains is to show that any element of $\ker \rho^*$
lies in the image of $\sigma_*$.
By Remark \ref{equivzero}
any element in the kernel
of $\rho^*$ must have a representative of the form $(f,\alpha)$
where, if we let
$g\colon S\to U$ denote the restriction of
$f\colon T\to X$,
there exist projective
morphisms $p_1,p_2\colon V\to S$ for some $V$ such that
$g\smallcirc p_1$ is 2-isomorphic to $g\smallcirc p_2$,
and
$\beta_i\in \widehat A^{p_i}_*V$ ($i=1,2$) such that
$\iota_{p_1}(\beta_1)=\iota_{p_2}(\beta_2)$
and such that
\begin{equation}
\label{cy}
\rho^{\prime{*}}\alpha=p_{2{*}}\beta_2-p_{1{*}}\beta_1
\end{equation}
in $\widehat A_*S$ (where $\rho'$ denotes the inclusion map $S\to T$).

Modifying $T$, we may assume the bundles on $S$ on whose
pullbacks $\beta_i$ are defined are in fact restrictions
of bundles $F_i$ on $T$.
Moreover, we may assume (\ref{cy}) holds as a rational equivalence
in a bundle which is the restriction of a bundle on $T$
(if there are surjections $G\to F_1|_S$ and $G\to F_2|_S$
then we form $0\to K\to G\oplus G\to (F_1\oplus F_2)|_S\to 0$ and
apply Corollary \ref{xtendxtensioncor}).
A choice of 2-isomorphism from $g\smallcirc p_1$ to $g\smallcirc p_2$
determines a morphism $V\to S\times_US$, which must be projective,
and by an application of Corollary \ref{xtendproj} we may assume we have a
fiber diagram
$$
\xymatrix{
Q \ar[r]^{\sigma''} \ar@<-2pt>[d]\ar@<2pt>[d]
\ar@/_1pc/[dd]_r &
W \ar@<-2pt>[d]_{q_1}\ar@<2pt>[d]^{q_2} &
V  \ar@<-2pt>[d]_{p_1}\ar@<2pt>[d]^{p_2} \ar[l] \\
R \ar[r]^{\sigma'} \ar[d]^h & T \ar[d]^f & S \ar[l]_{\rho'}\ar[d]^g \\
Y \ar[r]^\sigma & X & U \ar[l]_\rho
}
$$
After a modification of $W$ we may assume the equality
$\iota_{p_1}(\beta_1)=\iota_{p_2}(\beta_2)$ holds as a rational
equivalence in a bundle which is the restriction of a bundle on $W$.

We may find $\gamma_i\in \widehat A^{q_i}_* W$, represented by
a cycle in the bundle $q_i^*F_i$, which restricts to $\beta_i$,
for $i=1,2$.
By the excision axiom for na\"\i{}ve Chow groups,
$\iota_{q_2}(\gamma_2)-\iota_{q_1}(\gamma_1)=\sigma''_*\varepsilon$
for some $\varepsilon\in \widehat A_*^{\sigma''}Q$, and then
$\alpha=q_{2{*}}\gamma_2-q_{1{*}}\gamma_1+\sigma'_*\delta$
for some $\delta\in \widehat A^{\sigma'}_*R$.
Now, by Remark \ref{almostequiv},
$$(f,\alpha)=\sigma_*(r,\varepsilon)+\sigma_*(h,\delta). \qed $$
\renewcommand{\qed}{}\end{proof}

\section{Equivalence on bundles}
\label{vb}

\begin{conv}
If $E$ is a vector bundle on a stack $X$,
then $P(E)$ denotes the projectivization of $E$
in the sense of \cite{f},
i.e., $P(E)={\mathbf P}(\Sym^{\smallbullet}\cE^*)$,
where $\cE$ is the sheaf of sections of $E$.
If $E$ and $F$ are vector bundles on a stack $X$,
then $E\times_XF$ is also denoted using the Whitney sum notation
$E\oplus F$.
\end{conv}

Natural constructions in intersection theory produce cycles in
vector bundles.
An intersection theory requires, at the very least, a means of
intersecting with the zero section of a vector bundle.
This we develop by first examining trivial bundles (section \ref{trivlb}),
then developing a top Chern class operation (section \ref{topchernclass}),
followed by a construction which deals with general cycles in
vector bundles by compactifying and collapsing to the zero section
(sections \ref{collapsing} and \ref{izero}).
After discussing affine bundles (section \ref{affinebundles}),
we follow with
some additional topics which form a standard part of intersection
theory:
Segre classes, Chern classes, and the projective bundle theorem
(section \ref{segrechern}).

On a scheme $X$, with a vector bundle $\pi\colon E\to X$ of rank $e$
and associated projective bundle $p\colon P(E)\to X$,
the maps $\alpha\mapsto \pi^*\alpha$ and
$(\alpha_i)\mapsto \sum_{i=0}^{e-1}
c_1(\cO_E(1))^i\smallcap p^*\alpha_i$
are shown to be surjective using a noetherian
induction argument which reduces the claim to the case of
a trivial bundle; demonstrating injectivity takes more work
and is proved first for projectivized bundles, and then
for vector bundles (\cite[\S3]{f}).
The key fact that is used is that vector bundles on
schemes are Zariski locally trivial.
For stacks, the order of the steps is different: we first deduce
injectivity of $\pi^*$
using the explicit formula for $(\pi^*)^{-1}$ \cite[Proposition 3.3]{f}.
Surjectivity is harder, and uses a particular equivalence of
cycles on bundles (Proposition \ref{bundler}).
The projective bundle theorem follows using
the splitting principle (there is a pullback
which transforms an arbitrary bundle into a filtered bundle)
in a strong form (upon further pullback the filtration admits a splitting).

\subsection{Trivial bundles}
\label{trivlb}

Let us start with a fact which is familiar
from intersection theory on schemes.

\begin{pr}
\label{naivetriv}
Let $X$ be a stack, and $Z\subset X\times\A^1$ an integral
closed substack.
Let $\pi\colon X\times\A^1\to X$ be projection, with zero section $s$.
Then $[Z]$ is rationally equivalent to $\pi^*s^*[Z]$
on $X\times\A^1$.
\end{pr}

\begin{proof}
If $Z$ is contained in the fiber over zero then the
claim is clear, so assume the contrary.
Let $Z'$ be the closure in $X\times\A^1\times\bP^1$ of
the substack $(1_X\times m)^{-1}(Z)$ of
$X\times\A^1\times\A^1$, where
$m$ is multiplication, $m(x,u)=ux$.
If we think of $u$ as a time variable, then the
fiber at time 1 is a copy of $Z$, and as $u$ approaches 0,
the cycle $Z$ is stretched towards infinity, until at time $u=0$,
the fiber is identical to $\pi^{-1}(s^{-1}(Z))$.
Since $Z'\to \bP^1$ is flat away from infinity, the fibers over 0
and over 1 are rationally equivalent in $Z'$.
This rational equivalence pushes forward
to a rational equivalence between $\pi^*[s^{-1}(Z)]$ and $[Z]$ on
$X\times\A^1$.
\end{proof}

\relax From this, it is straightforward to deduce the corresponding result
for $A_*$.

\begin{pr}
\label{trivlbpr}
Let $X$ be a stack, and let $\pi\colon X\times\A^1\to X$ be
the projection map with zero section $s$.
Then $\pi^*\colon A_kX\to A_{k+1}X\times\A^1$ is an
isomorphism, with inverse $s^*$.
\end{pr}

\begin{proof}
Clearly we have $s^*\smallcirc \pi^* = {\rm id}$.
It remains, therefore, to show that $\pi^*$ is surjective.
By Proposition \ref{exci}, it suffices to show that any cycle class
on $X\times\bP^1$ is represented by a cycle whose restriction to
$X\times\A^1$ lies in the image of $\pi^*$.  So, let
$(f,\alpha)$ be a cycle on $X\times\bP^1$,
with $f\colon T\to X\times\bP^1$ projective.
We may assume $T$ is integral and $\alpha$ is represented by
an integral closed substack $Z$ of the total space of
a vector bundle $E\to T$.
The cases where the composite $Z\to\bP^1$ has image equal to
$\{0\}$ or $\{\infty\}$ are easily dealt with, so we assume
the contrary.

Let $Y=X\times\bP^1$.
Let $P$ be the closure in $(Y\times\bP^1)\times_XY$ of the
graph $\Gamma$ of the morphism
$1_X\times m\colon X\times\A^1\times\A^1\to X\times\A^1$
(where $m$ is as in Proposition \ref{naivetriv}).
We have projections $\psi\colon P\to Y\times\bP^1$
and $\varphi\colon P\to Y$.
For $t$ in the base field, we may form the fiber diagram
$$\xymatrix{
{\widetilde E_t}\ar[r]\ar[d] &
{\widetilde T_t}\ar[r]^{f_t} \ar[d] &
P_t \ar[r]^(.35){\psi_t} \ar[d]^{i_t} & Y\times\{t\}\ar[d] \\
{\widetilde E}\ar[r]\ar[d] & {\widetilde T}\ar[r] \ar[d] &
P \ar[r]^(.35)\psi \ar[d]^\varphi & Y\times\bP^1 \\
E \ar[r] & T \ar[r]^f & Y
}$$

Although
the map $\varphi$ is not flat, its restriction to $\Gamma$ is flat.
If $[(\varphi|_\Gamma)^{-1}(Z)]=\sum a_i [W_i]$, we let
$\tilde \alpha=\sum a_i [\overline{W_i}]$
and $\tilde \alpha_t=\sum a_i [\overline{W_i}\times_{\bP^1}\{t\}]$.
Then
$(\psi_0\smallcirc f_0,\tilde \alpha_0)=
(\psi_1\smallcirc f_1,\tilde \alpha_1)$
in $A_*Y$.
We now recognize that
$\psi_1$ and $\varphi\smallcirc i_1$ are
isomorphisms, and that if we form the fiber square
$$\xymatrix{
P'_0 \ar[r]^(.3){\psi'_0} \ar[d] & X\times\A^1\times\{0\}\ar[d] \\
P_0 \ar[r]^(.35){\psi_0} & Y\times\{0\}
}$$
then ${\psi'_0}$ is an isomorphism,
and the composite $X\times\A^1\eqto P'_0\to P_0\to P\to Y$
is equal to projection $X\times\A^1\to X$ followed by
inclusion as the zero section $X\to Y$.
Thus $(\psi_0\smallcirc f_0,\tilde \alpha_0)$ restricts over $X\times\A^1$ to
a pullback from $X$.
\end{proof}

\subsection{The top Chern class operation}
\label{topchernclass}

\begin{defn}
\begin{romanilist}
\firstitem Let $X$ be a connected stack,
and let $U\to X$ be a vector bundle of rank $r$.  Then we have a map
$$\widehat A_jX\to \widehat A_{j-r}X$$
which, for any vector bundle $E\to X$, with $\rk E=s$,
sends $\alpha\in A^\circ_{j+s}E$ to
$s_*\alpha\in A^\circ_{j+s}E\oplus U$,
where $s$ is the zero section of $E\oplus U\to E$.
We denote the image class in $\widehat A_{j-r}X$ by $\alpha^U$.
For $\alpha\in\widehat A_*X$ where $X$ is not connected,
$\alpha^U$ denotes the sum of $(\alpha|_{X_i})^U$ over
connected components $X_i$ of $X$.
\item Let $Y$ be a stack, and let $U\to Y$ be a vector
bundle of rank $r$.
We define the {\em top Chern class operation}
$$c_{\rm top}(U)\smallcap{-}\colon A_jY\to A_{j-r}Y$$
by $(f,\alpha)\mapsto(f,\alpha^{f^*U})$.
That this map respects equivalence is clear.
\end{romanilist}
\end{defn}

\begin{pr}
\label{cycleinbundleprop}
Let $Y$ be a stack, $\pi\colon U\to Y$ a vector bundle,
$f\colon T\to Y$ a projective morphism, $E\to T$ a vector bundle,
and suppose $\gamma\in A^\circ_*E\oplus f^*U$.
If we form the fiber diagram
$$
\xymatrix{
U \ar[d]_\pi & f^*U \ar[l]_{f'} \ar[d] & E\oplus f^*U \ar[l] \ar[d] \\
Y & T \ar[l]_f & E \ar[l]
}
$$
then $\pi^*(f,\gamma)=(f',\gamma)$,
where on the right we view $\gamma$ as a cycle in the
total space of the bundle $E\oplus f^*U\to f^*U$.
\end{pr}

\begin{proof}
If $\tau$ denotes projection onto the first two factors
$E\oplus f^*U\oplus f^*U\to E\oplus f^*U$
then
$\pi^*(f,\gamma)=(f',\tau^*\gamma)$.
But now we may let $\upsilon$ denote projection onto the
first and third factors
$E\oplus f^*U\oplus f^*U\to E\oplus f^*U$;
by Remark \ref{twovbsurj},
$\tau^*\gamma=\upsilon^*\gamma$ in
$A^\circ_*E\oplus f^*U\oplus f^*U$,
and we are done
($\upsilon^*\gamma\in A^\circ_*E\oplus f^*U\oplus f^*U$
and $\gamma\in A^\circ_*E\oplus f^*U$
determine the same element of $\widehat A_*f^*U$).
\end{proof}

\begin{cor}
\label{topchernprop}
Let $Y$ be a stack, and let $\pi\colon U\to Y$
be a vector bundle with zero section $s$.
Then
$\pi^*c_{\rm top}(U)\smallcap\alpha=s_*\alpha$
for all $\alpha\in A_*Y$.
\end{cor}

The top Chern class operation obeys the expected properties.
\begin{pr}
\label{expectedprop}
\begin{romanilist}
\firstitem If $E$ and $F$ are vector bundles on a stack $X$,
then we have $c_{\rm top}(E)\smallcap (c_{\rm top}(F)\smallcap \alpha) =
c_{\rm top}(F)\smallcap (c_{\rm top}(E)\smallcap \alpha)$
for all $\alpha\in A_*X$.
\item If $f\colon Y\to X$ is projective and $E$ is a vector
bundle on $X$, then
$f_*(c_{\rm top}(f^*E)\smallcap \alpha) = c_{\rm top}(E)\smallcap f_*\alpha$
for all $\alpha\in A_*X$.
\item If $f\colon Y\to X$ is flat of locally constant relative dimension
and $E$ is a vector bundle on $X$, then
$c_{\rm top}(f^*E)\smallcap f^*\alpha = f^*(c_{\rm top}(E)\smallcap \alpha)$
for all $\alpha\in A_*X$.
\item If $E$ is a vector bundle on a stack $X$ with a nowhere
vanishing section $s$, then $c_{\rm top}(E)\smallcap \alpha = 0$
for all $\alpha\in A_*X$.
\item If $L_1$ and $L_2$ are line bundles on a stack $X$, then
$c_{\rm top}(L_1)\smallcap\alpha + c_{\rm top}(L_2)\smallcap\alpha =
c_{\rm top}(L_1\otimes L_2)\smallcap\alpha$ for all $\alpha\in A_*X$.
\item If $L$ is a line bundle on a stack $X$, then
$c_{\rm top}(L)\smallcap\alpha + c_{\rm top}(L^\vee)\smallcap\alpha = 0$
for all $\alpha\in A_*X$.
\end{romanilist}
\end{pr}

\begin{proof}
Routine,
for instance, the content of (v) is an explicit rational
equivalence on $L_1\oplus L_2\oplus (L_1\otimes L_2)$,
given locally by $(x,y,z)\mapsto xyz^{-1}$.
\end{proof}

\subsection{Collapsing cycles to the zero section}
\label{collapsing}

We present here a construction which is the reverse process of
that of Section \ref{trivlb}.
Instead of stretching a cycle to infinity,
we collapse a general cycle on a projectivized bundle $P(U\oplus 1)$
to the zero section, obtaining, away from the zero section, the
pullback of the intersection with the boundary $P(U)$.

We first need some basic geometry of bundles and
their projectivizations.
Let $Y$ be a stack, with vector bundle $U$ of rank $r$, and
consider the diagram as above:
$$\xymatrix{
P(U)\ar[r]^(.42)i\ar[dr]_\tau & P(U\oplus 1)\ar[d]_(.45)\sigma &
U\ar[l]_(.35)\rho \ar[dl]_\pi \\
& Y \ar@/_1pc/[ur]_s
}$$
Let $s_0:=\rho\smallcirc s$ be the zero section of $P(U\oplus 1)\to Y$,
and let $\sigma_0:=s_0\smallcirc\sigma$ be the map collapsing
all of $P(U\oplus 1)$ onto the zero section.
The map which projects $P(U\oplus 1)$ onto $P(U)$ is only
a rational map, but
becomes regular after blowing up
the zero section.
Denote the resulting projection morphism
$\Bl_{s_0(Y)}P(U\oplus 1)\to P(U)$ by $\eta$, and let $\xi$ be
the blow-down.
The maps $\eta$ and $\xi$ determine a closed immersion
\begin{equation}
\label{emb}
\Bl_{s_0(Y)}P(U\oplus 1)\to P(U)\times_Y P(U\oplus 1)
\end{equation}

\begin{pr}
\label{cyclecollapse}
With notation as above,
let $f\colon T\to P(U\oplus 1)$ be projective and
let $(f,\alpha)$ be a cycle on $P(U\oplus 1)$ represented by
$[Z]\in A^\circ_*E$ for some vector bundle $E$ on $T$,
with $Z$ integral and not supported over
$P(U)$.
Let $f'\colon T'\to P(U)$ be the restriction of $f$ and let
$\alpha'\in \widehat A_*T'$
be given by $[Z\times_{P(U\oplus 1)}P(U)]\in A^\circ_*E|_{T'}$.
Then
\begin{equation}
\label{collapsecycle}
(f,\alpha)=\xi_*\eta^*(f',\alpha')+\sigma_{0{*}}(f,\alpha)
\end{equation}
in $A_*P(U\oplus 1)$.
\end{pr}

\begin{proof}
Let $S=\Bl_{s_0(Y)\times\{0\}\cup P(U)\times\{\infty\}}
P(U\oplus 1)\times\bP^1$.
We have an inclusion
$$S\hookrightarrow (P(U\oplus 1)\times\bP^1)\times_Y P(U\oplus 1)$$
whose image is given locally in coordinates as
\begin{multline*}
\{\,((x_1:\ldots:x_r:t),(u:v))\times(X_1:\ldots:X_r:T)\mid \\
(X_1,\ldots,X_r)\propto (x_1,\ldots,x_r){\rm\ and\ }
(X_1,\ldots,X_r,T)\propto (ux_1,\ldots,ux_r,tv)\,\}.
\end{multline*}
Let $\psi\colon S\to P(U\oplus 1)\times\bP^1$
be the first projection; $\psi$ is the blow-down map.
Let $\varphi$ be the second projection map.

The crucial properties of these maps
are that away from $\infty=(0:1)$ in $\bP^1$,
the map $\varphi$ is flat,
and that the fiber of $S$ over $0=(1:0)\in\bP^1$ is
isomorphic to $(\Bl_{s_0(Y)}P(U\oplus 1))\cup P(U\oplus 1)$,
where the first component is blown down by $\psi$
onto $P(U\oplus 1)\times\{0\}$ and mapped
by $\varphi$ via the projection map $\eta$ above,
and the second component is collapsed to the zero section
by $\psi$ and mapped isomorphically by $\varphi$.
As in section \ref{trivlb} we think of $(u:v)$ as
a time variable, and at time
$1=(1:1)\in\bP^1$, both $\psi$ and $\varphi$ restrict to
isomorphisms with $P(U\oplus 1)$.

The idea, expressed in the following diagram, is that
the operation of pulling back a cycle on $P(U\oplus 1)$
via $\varphi$ and projecting by $\psi$ produces
the pullback to $P(U\oplus 1)$ of the intersection of the cycle
with the boundary $P(U)\subset P(U\oplus 1)$, plus a
copy of the projection of the cycle to the zero section.
We consider the fiber diagram
$$\xymatrix{
{\widetilde E_t}\ar[r]\ar[d] &
{\widetilde T_t}\ar[r]^{f_t} \ar[d] &
S_t \ar[r]^(.35){\psi_t} \ar[d]^{i_t} & P(U\oplus 1)\times\{t\}\ar[d] \\
{\widetilde E}\ar[r]\ar[d] & {\widetilde T}\ar[r] \ar[d] &
S \ar[r]^(.35)\psi \ar[d]^\varphi & P(U\oplus 1)\times\bP^1 \\
E \ar[r] & T \ar[r]^(.4)f & P(U\oplus 1)
}$$
and then proceeding just as in the proof of Proposition \ref{trivlbpr},
we deduce
\begin{equation}
\label{collap}
(\psi_1\smallcirc f_1, [\widetilde Z_1]\in A^\circ_* \widetilde E_1)=
(\psi_0\smallcirc f_0, [\widetilde Z_0]\in A^\circ_* \widetilde E_0)
\end{equation}
in $A_*P(U\oplus 1)$, where
$\widetilde Z_t=\varphi^{-1}(Z)\times_{\bP^1}\{t\}$.
The left- and right-hand sides of (\ref{collap}) are
equal to the left- and right-hand sides
of (\ref{collapsecycle}), respectively.
\end{proof}

\subsection{Intersecting cycles with the zero section}
\label{izero}

We finally give a proof
that the pullback map $\pi\colon U\to X$ for a vector
bundle induces an isomorphism on $A_*$.
This is achieved by giving an explicit inverse to $\pi^*$.
We continue to use the notation of
the previous section:
$\pi\colon U\to Y$ is a vector bundle,
with morphisms $\sigma$, $\tau$, etc.\ as before.

\begin{pr}
\label{bundler}
Let $R$ be the universal quotient bundle of $\tau^*U$
on $P(U)$.
Then for any $\alpha\in A_*P(U)$,
$\pi^*\tau_* c_{\rm top}(R)\smallcap\alpha=
\rho^*\xi_*\eta^*\alpha$ in $A_*U$.
\end{pr}

\begin{proof}
Let $\alpha$ be represented by $(f,[Z])$, where
$f\colon T\to P(U)$ is projective and $Z$ is
an integral closed substack of $E$ for some
vector bundle $E\to T$.
Form the fiber diagram
$$\xymatrix{
U \ar[d]^\pi & Q \ar[l]_\upsilon \ar[d] &
S \ar[l]_g \ar[d]^{\pi'} & F \ar[l] \ar[d]^{\pi''} \\
Y & P(U) \ar[l]_\tau & T \ar[l]_(.4)f & E \ar[l]
}$$

By the exact sequence of bundles on
$P(U)$
$$0\to \cO_U(-1)\to \tau^*U\to R\to 0,$$
we have
$$c_{\rm top}(R)\smallcap \alpha =
(f,[Z\times_Tf^*\cO_U(-1)]\in A^\circ_*F).$$
Then
$$\tau_*c_{\rm top}(R)\smallcap \alpha =
(\tau\smallcirc f, [Z\times_Tf^*\cO_U(-1)]\in A^\circ_*F).$$
Now by Proposition \ref{cycleinbundleprop} we have
$$\pi^*\tau_*c_{\rm top}(R)\smallcap \alpha =
(\upsilon\smallcirc g, [Z\times_Tf^*\cO_U(-1)]\in A^\circ_*F),$$
and this last expression we recognize as being a representative for
$\rho^*\xi_*\eta^*\alpha$.
\end{proof}

\begin{pr}
Let $Q$ be the universal quotient bundle of
$\sigma^*(U\oplus 1)$ on $P(U\oplus 1)$.
Then
\begin{romanlist}
\item for any $\alpha\in A_*P(U)$,
$c_{\rm top}(Q)\smallcap i_*\alpha=0$;
\item for any $\alpha\in A_*P(U\oplus 1)$,
$\pi^*\sigma_*c_{\rm top}(Q)\smallcap \alpha = \rho^*\alpha$.
\end{romanlist}
\end{pr}

\begin{proof}
Let $f\colon T\to P(U\oplus 1)$ be projective,
let $E\to T$ be a vector bundle,
let $Z$ be an integral closed substack of $E$,
and let $\alpha\in\widehat A_*T$ be given by
$[Z]\in A^\circ_*E$.

By the exact sequence of bundles on $P(U\oplus 1)$
$$0\to \cO_{U\oplus 1}(-1)\to \sigma^*(U\oplus 1)\to Q\to 0,$$
$\alpha^{f^*Q}\in\widehat A_*T$ is
represented by
$[Z\times_T f^*\cO_{U\oplus 1}(-1)]\in A^\circ_*E\oplus f^*\sigma^*U\oplus 1$.
If we consider $E\oplus f^*\sigma^*U\oplus 1\to E\oplus f^*\sigma^*U$ as
a trivial $\A^1$-bundle with zero section $t$
then by Proposition \ref{naivetriv},
$\alpha^{f^*Q}$ is also represented
by $t^*[Z\times_T f^*\cO_{U\oplus 1}(-1)]\in A^\circ_*E\oplus f^*\sigma^*U$.
The intersection of $\cO_{U\oplus 1}(-1)$ with the
zero section of $\sigma^*U\oplus 1\to\sigma^*U$ is
the union of $\cO_{U\oplus 1}(-1)|_{P(U)}(=\cO_U(-1))$ and
the zero section of $\cO_{U\oplus 1}(-1)$.
Statement (i) follows since $t^*$ is defined to be
zero on cycles contained in the support of the divisor.

In the remaining case, that is, when the projection of $Z$
to $P(U\oplus 1)$ does not factor through $P(U)$,
let $(f',\alpha':=[Z']\in A^\circ_*E|_{T'})$
be determined from $(f,\alpha:=[Z]\in A^\circ_*E)$
as in the statement of Proposition \ref{cyclecollapse}.
Then
\begin{equation}
\label{tstar}
t^*[Z\times_T f^*\cO_{U\oplus 1}(-1)]=
[Z'\times_{T'} f^{\prime{*}}\cO_U(-1)]+
\tilde s_*[Z],
\end{equation}
in $A^\circ_* E\oplus f^*\sigma^*U$, where $\tilde s$ is the zero
section of $E\oplus f^*\sigma^*U\to E$.
Now, if we let $R$ be the universal quotient bundle
of $\tau^*U$ on $P(U)$, then (\ref{tstar}) implies the formula
\begin{equation}
\label{vbequiv}
c_{\rm top}(Q)\smallcap (f,\alpha)=
i_*c_{\rm top}(R)\smallcap (f',\alpha') +
c_{\rm top}(\sigma^*U)\smallcap (f,\alpha).
\end{equation}
By Proposition \ref{bundler} and the projection formula coupled with
Corollary \ref{topchernprop}, we now have
$$\pi^*\sigma_* c_{\rm top}(Q)\smallcap (f,\alpha) =
\rho^*\xi_*\eta^*(f',\alpha')+s_*\sigma_*(f,\alpha),$$
and the desired statement follows by Proposition \ref{cyclecollapse}.
\end{proof}

\begin{pr}
\label{tctt}
Let $\pi\colon U\to Y$ be a vector bundle,
with projectivization $\tau\colon P(U)\to Y$.
Let $R$ be the universal quotient bundle of $\tau^*U$.
Then $$\tau_*c_{\rm top}(R)\smallcap \tau^*\alpha=\alpha$$
for all $\alpha\in A_*Y$.
\end{pr}

\begin{proof}
If $\alpha$ is represented by $(f,[Z]\in A^\circ_*E)$
where $f\colon T\to Y$ is projective and
$E\to T$ is a vector bundle,
and we form the fiber square
$$\xymatrix{
{\widetilde T} \ar[r]^{\tilde\tau}\ar[d]^{\tilde f} & T \ar[d]^f \\
P(U)\ar[r]^(.6)\tau & Y
}$$
then $c_{\rm top}(R)\smallcap \tau^*\alpha$ is
represented by
$[Z\times_T\tilde f^*\cO_U(-1)]\in A^\circ_*(\tilde\tau^*E\oplus
\tilde\tau^*f^*U)$.
Since the bundle on which this cycle lives is
the pullback of a bundle on $T$, we may compute
$\tau_*c_{\rm top}(R)\smallcap \tau^*\alpha$
by pushing forward the actual cycle.
The push-forward map, restricted to
$Z\times_{\widetilde T}\tilde f^*\cO_U(-1)$,
is a birational map,
and hence $\tau_*c_{\rm top}(R)\smallcap \tau^*\alpha$
is represented by
$[Z\times_Tf^*U]\in A^\circ_*(E\oplus f^*U)$, i.e., is
equal to $\alpha$.
\end{proof}

\begin{cor}
\label{pistariso}
The map $\alpha\mapsto \sigma_* c_{\rm top}(Q)\smallcap \bar\alpha$,
where $\bar \alpha\in A_*P(U\oplus 1)$ is some cycle class
that restricts to $\alpha$ (guaranteed to exist by excision),
is independent of the choice
of $\bar\alpha$
and gives an isomorphism $A_*U\to A_*Y$,
and this map is inverse to the map $\pi^*$.
\end{cor}

\subsection{Affine bundles}
\label{affinebundles}

The structure group for affine $n$-plane bundles
is the subgroup $\Aff(n)$ of $GL(n+1)$ of matrices which have
$n$ zeros followed by a 1 along the bottom row.
Let $X$ be a stack;
the map $\Aff(n)\to GL(n+1)$ associates to every locally trivial
(for the smooth topology) affine $n$-plane bundle
$\varphi\colon B\to X$
a vector bundle $E\to X$ of rank $(n+1)$, together with a
surjection of vector bundles $\tau\colon E\to \cO_X$.
If $F$ denotes the kernel to $\tau$ then we have a closed immersion
$P(F)\to P(E)$ with complement isomorphic to $B$.
The affine bundle $\varphi\colon B\to X$ has the property that
the associated surjection $\tau\colon E\to \cO_X$
admits a splitting after pullback via $\varphi$.

Since a general affine bundle has no zero section,
the formalism of section \ref{collapsing} does not apply,
although we may still use the argument of section \ref{izero}
to deduce that $\varphi^*$ is a split monomorphism.
That $\varphi^*$ is a split monomorphism
is one ingredient in the proof of the projective bundle theorem
(section \ref{segrechern}).
Surjectivity of $\varphi^*$ comes only later as a corollary to the
projective bundle theorem.

\begin{lm}
\label{istarzero}
Let $X$ be a scheme, and let $0\to F\stackrel{i}\to E\to \cO_X\to 0$ be
an exact sequence of vector bundles on $X$.
Then $i_*\alpha=0$ for all $\alpha\in A_*F$.
\end{lm}

\begin{proof}
Denote by $\pi_E$ and $\pi_F$ the respective projections from $E$ and $F$.
The Gysin map $i^*$ (from section \ref{basic}) satisfies
$i^*\smallcirc\pi_E^*=\pi_F^*$, so
by Corollary \ref{pistariso}, $i^*$ is an isomorphism.
By the definition of $i^*$ we have $i^*i_*\alpha=0$, so we conclude
that $i_*\alpha=0$.
\end{proof}

\begin{pr}
Let $\varphi\colon B\to X$ be an affine $n$-plane bundle with associated
vector bundle $E$ and exact sequence
$$0\to F\to E\to \cO_X\to 0$$
of vector bundles on $X$.
Denote by $Q$ the universal quotient bundle of the pullback of $E$
via $\sigma\colon P(E)\to X$.
Let $i$ denote the map $P(F)\to P(E)$.
Then $c_{\rm top}(Q)\smallcap i_*\alpha=0$ for all $\alpha\in A_*P(F)$.
\end{pr}

\begin{proof}
Consider the pullback sequence
$0\to \sigma^*F\stackrel{i'}\to\sigma^*E\to \cO_{P(E)}\to 0$.
With the diagram
$$\xymatrix{
{\cO_E(-1)}\ar[r]^(.55)j\ar[dr]_\mu &
{\sigma^*E} \ar[r]\ar[d]_(.45)\nu & Q \\
& P(E)
}$$
we have
$\nu^*(c_{\rm top}(Q)\smallcap i_*\alpha)=j_*\mu^*i_*\alpha=
i'_*k_*\psi^*\alpha$,
where $\psi$ denote the projection $\cO_F(-1)\to P(F)$ and
$k$ denote the inclusion of $\cO_F(-1)$ in $\sigma^*F$.
But $i'_*k_*\psi^*\alpha=0$ by Lemma \ref{istarzero},
and $\nu^*$ is an isomorphism by Corollary \ref{pistariso},
so we are done.
\end{proof}

Now Proposition \ref{tctt} gives us
\begin{cor}
\label{pistarsplitmono}
With notation as above, $\varphi^*$ is a split monomorphism.
A splitting is the map sending
$\alpha\in A_*B$ to
$\sigma_*(c_{\rm top}(Q)\smallcap \bar\alpha)$,
where $\bar\alpha$ is any cycle class on $P(E)$ which restricts to $\alpha$.
\end{cor}

\subsection{Segre and Chern classes and the projective bundle theorem}
\label{segrechern}

It is now routine to define general Segre classes and Chern classes
of vector bundles.
The usual projective bundle theorem will appear as a consequence.

\begin{defn}
Let $\pi\colon E\to X$ be a vector bundle on a stack $X$.
The $i^{\rm th}$ Segre class operation
$s_i(E)\smallcap\colon A_kX\to A_{k-i}X$ is
defined by the formula
$$s_i(E)\smallcap\alpha =
p_*( c_{\rm top}(\cO_E(1))^{\rk E - 1 + i}\smallcap p^*\alpha),$$
where $p$ denotes the projection $P(E)\to X$.
\end{defn}

\begin{pr}
\label{segreprop}
\begin{romanslist}
\firstitem Let $E$ be a vector bundle on a stack $X$.
Then for all $\alpha\in A_*X$, we have
$s_i(E)\smallcap\alpha = 0$ for $i<0$ and
$s_0(E)\smallcap\alpha = \alpha$.
\item If $E$ and $F$ are vector bundles on a stack $X$
and $\alpha\in A_*X$, then for all $i$ and $j$,
$s_i(E)\smallcap (s_j(F)\smallcap\alpha) =
s_j(F)\smallcap (s_i(E)\smallcap \alpha)$.
\item If $f\colon X'\to X$ is projective, $E$ is a vector bundle
on $X$, and $\alpha\in A_*X'$, then, for all $i$,
$f_*(s_i(f^*E)\smallcap\alpha)=s_i(E)\smallcap f_*\alpha$.
\item If $f\colon X'\to X$ is flat of locally constant relative dimension,
$E$ is a vector bundle on $X$, and $\alpha\in A_*X$, then, for all $i$,
$s_i(f^*E)\smallcap f^*\alpha = f^*(s_i(E)\smallcap \alpha)$.
\item If $E$ is a line bundle on $X$ and
$\alpha\in A_*X$, then
$s_1(E)\smallcap\alpha = -c_{\rm top}(E)\smallcap \alpha$.
\end{romanslist}
\end{pr}

\begin{proof}
Only part (i) is nontrivial.
By Corollary \ref{pistariso}, it is enough to verify that the identities
of part (i) hold
after pullback to vector bundles.
Consider the fiber diagram
$$\xymatrix{
P(\pi_1^*E)\ar[r]^(.6){q_1} \ar[d]^{\sigma_1} & E^\vee \ar[d]^{\pi_1} \\
P(E) \ar[r]^(.55)p & X
}$$
and the sequence of bundles
$$0\to S_1\stackrel{\iota_1}\to P(\pi_1^*E)\to \cO_E(1)\to 0$$
on $P(\pi_1^*E)$.
Inductively, we define $\sigma_k$, $\pi_k$, and $q_k$ by the
fiber diagram
$$\xymatrix{
P(\pi_k^*\pi_{k-1}^*\cdots\pi_1^*E)\ar[r]^(.68){q_k} \ar[d]^{\sigma_k} &
{E^\vee}^{\oplus k} \ar[d]^{\pi_k} \\
P(\pi_{k-1}^*\cdots\pi_1^*E)\ar[r]^(.65){q_{k-1}} & {E^\vee}^{\oplus k-1}
}$$
with
$$\xymatrix{
S_k\ar[r]^(.35){\iota_k}\ar[dr]_{\tau_k} & P(\pi_k^*\cdots\pi_1^*E)
\ar[r]\ar[d]_(.45){\sigma_k} & \sigma_{k-1}^*\cdots\sigma_1^*\cO(1) \\
& P(\pi_{k-1}^*\cdots\pi_1^*E)
}$$
We have
$$\pi_k^*\cdots\pi_1^*(c_1 \cO(1)^k\smallcap p^*\alpha) =
q_{k{*}}(\iota_{k{*}}\tau_k^*)\cdots(\iota_{1{*}}\tau_1^*) p^*\alpha.$$
Set $e=\rk E - 1$;
it suffices to show that
\begin{equation}
\label{segreid}
q_{k{*}}(\iota_{k{*}}\tau_k^*)\cdots(\iota_{1{*}}\tau_1^*) p^*\alpha =
\begin{cases}
0 & \text{if $k<e$,} \\
\pi_k^*\cdots\pi_1^*\alpha & \text{if $k=e$.}
\end{cases}
\end{equation}

What we do is show that the identity (\ref{segreid}) holds for $Z_*$.
This is a local computation, so we may assume $E$ is the trivial bundle
of rank $e+1$.
The cycle we are pushing forward along
$X\times\bP^e\times (\A^{e+1})^k\to X\times (\A^{e+1})^k$
has fiber equal to the set of
points in $\bP^e$ satisfying $k$ linear conditions, which are
generically independent.
Thus the generic fiber is positive-dimensional for $k<e$ and
has degree 1 when $k=e$.
\end{proof}

\begin{pr}
Let $\pi\colon E\to X$ be a vector bundle on a stack $X$
with $\rk E = e+1$.
Let $\sigma$ denote the projection $P(E)\to E$.
Then the map $\theta_E\colon A_*X^{e+1}\to A_*P(E)$
given by $(\alpha_0,\ldots,\alpha_e)\mapsto
\sum_{i=0}^e c_1(\cO_E(1))^i\smallcap \sigma^*\alpha_i$
is an isomorphism.
\end{pr}

\begin{proof}
Injectivity is clear by Proposition \ref{segreprop}.
To demonstrate surjectivity of $\theta_E$, we induct on the
rank of $E$.
The case of a line bundle is trivial.
For the inductive step, we first of all consider
the special case where $E$ is filtered
as $E\to L\to 0$, for
some line bundle $L$ on $X$.
It suffices to demonstrate that $\theta_{E\otimes L^\vee}$ is
surjective
(we have $\nu\colon P(E)\stackrel{\raise-1.2pt\hbox{$\scriptstyle\sim$}}\to
P(E\otimes L^\vee)$ with $\nu^*\cO_{E\otimes L^\vee}(1)=
\cO_E(1)\otimes\sigma^*L$), so
we are reduced to the case where we have $E\to \cO_X\to 0$.
Then $E$ is associated to an affine bundle $\varphi\colon B\to X$,
and we have seen that $\varphi^*E\to\cO_B\to 0$ admits a splitting.
Writing $\varphi^*E\simeq F\oplus 1$ we have
$P(F)\to P(F\oplus 1)$ with complement $F$, so by the induction hypothesis and
the exact sequence
$$A_*P(F)\to A_*P(F\oplus 1)\to A_*F\to 0$$
we obtain that $\theta_{\varphi^*E}$ is surjective.
By Corollary \ref{pistarsplitmono} and Proposition \ref{expectedprop}
we conclude that $\theta_E$ is surjective.

We deduce the general case from the case above
by the standard splitting construction:
after pullback by $\tau\colon P(E^\vee)\to X$ we have
$\tau^*E\to \cO_{E^\vee}(1)\to 0$,
and as before (this time using Proposition \ref{tctt}), we deduce surjectivity
of $\theta_E$ from surjectivity of $\theta_{\tau^*E}$.
\end{proof}

We define Chern classes in terms of Segre classes using the universal
polynomials.
To deduce standard facts about Chern classes
(vanishing of $c_i(E)$ for $i>\rk E$,
projection formula, pullback formula,
and Whitney sum formula)
we use the splitting construction
to reduce to the case of direct sums of line bundles,
and then the formula
$c_1(L)=c_{\rm top}(L)$ (Proposition \ref{segreprop} (v)) implies the
desired statements.
By the same method, the formulas of
\cite[Remark 3.2.3]{f} (Chern classes of dual bundles,
tensor produces, etc.) hold as well for stacks; for instance we deduce
the formula
$$\zeta^r + c_1(p^*E)\zeta^{r-1} + \cdots +
c_r(p^*E)=0$$
which characterizes $A_*P(E)$, where $E$ is a vector bundle
of rank $r$ on a stack $X$ with projection
$p\colon P(E)\to X$ and with $\zeta=c_1(\cO_E(1))$
(\cite[Remark 3.2.4]{f}).

Let $\pi\colon E\to X$ be a vector bundle of rank $e$, with
zero section $s$.
The identity
$\pi^*c_{\rm top}(E)\smallcap \alpha = s_*\alpha$
(Corollary \ref{topchernprop})
which characterizes $c_{\rm top}(E)$ also holds for $c_e(E)$,
from which we deduce $c_{\rm top}(E)=c_e(E)$,
and so from now on we may use the two notations interchangeably.

Finally, the projective bundle theorem implies that the pullback
map on Chow groups induced by an affine bundle is an isomorphism.

\begin{cor}
\label{affinebundlecor}
Let $\varphi\colon B\to X$ be an affine bundle.
Then $\varphi^*\colon A_*X\to A_*B$ is an isomorphism.
\end{cor}

\begin{proof}
Let $E$ and $F$ be the associated vector bundles,
as in section \ref{affinebundles}.
Let $i$ denote the map $P(F)\to P(E)$,
and let $p$ and $q$ denote the projection to $X$ from
$P(E)$ and from $P(F)$, respectively.
We observe that $\cO_E(1)$ has a section vanishing precisely
on $P(F)$,
which implies
$c_1(\cO_E(1))\smallcap p^*\alpha=i_*q^*\alpha$
and hence
$c_1(\cO_E(1))^j\smallcap p^*\alpha=i_*c_1(\cO_F(1))^{j-1}\smallcap q^*\alpha$
for $j\ge 1$.
That $\varphi^*$ is an isomorphism now follows by the
projective bundle theorem plus the excision axiom.
\end{proof}

\begin{rem}
\label{affinebundlerem}
On a scheme, every affine bundle is locally trivial for
the Zariski topology,
so surjectivity of $\varphi^*$ follows by the elementary
argument of \cite[\S1.9]{f},
and then injectivity can be proved using the same argument
as for vector bundles (as remarked in \cite{semc}, this argument uses
only the existence of the compactification $P(E)$).
More advanced techniques, cf.\ \cite{gadv}, demonstrate homotopy invariance
for any flat bundle on a scheme base such that the fibers
over all points (closed or not) are isomorphic to affine spaces.
\end{rem}

\section{Elementary intersection theory}
\label{elementary}

\subsection{Fulton-MacPherson construction for local immersions}
\label{fmp}

The existence of Gysin maps for principal effective
Cartier divisors,
together with the homotopy property (Corollary \ref{pistariso}) and
excision axiom,
lets us apply the standard Fulton-MacPherson construction \cite[\S6]{f} to
produce Gysin maps for regular local immersions of Artin stacks.
We recall that a representable morphism $f\colon F\to G$ of
stacks is unramified if and only if, for some or equivalently
every smooth atlas $V\to G$, the pullback
$\tilde f\colon F\times_GV\to V$
fits into a commutative diagram
$$
\xymatrix{
S\ar[r]^{\tilde g}\ar[d] & T \ar[d] \\
F\times_GV\ar[r]^(.6){\tilde f} & V
}
$$
with $\tilde g$ a closed immersion and the vertical maps \'etale surjective
(\cite[Corollary 18.4.8]{egafour} plus the local nature of the
property of being unramified).
A representable morphism is called a {\em local immersion} if it is
unramified, in which case there is a well-defined normal cone,
given locally as the normal cone to $S\to T$.
A representable morphism is called a {\em regular local immersion}
if it is a local immersion and if moreover the map $\tilde g$
is a regular immersion (or equivalently, if the normal cone is
a vector bundle).

Given a local immersion $f\colon F\to G$, and given an
arbitrary morphism
$g\colon G'\to G$,
there corresponds a map
$A_*G'\to A_*(C_FG\times_FF')$.
If we set $F'=F\times_GG'$,
then $f'\colon F'\to G'$ is also a local immersion,
and we may form the deformation space
$M^\circ_{F'}G'\to\bP^1$,
which has general fiber $G'$ and
special fiber $s\colon C_{F'}G'\to M^\circ_{F'}G'$
(cf.~\cite[\S5.1]{f} for the case of a closed immersion,
and \cite{kr} for the case of a local immersion).
The map $A_*G\to A_*(C_FG\times_FF')$ is defined to be the composite
\begin{align*}
A_*G' & \to A_{*+1}G'\times(\bP^1\setminus\{0\}) \\
& \eqto A_{*+1}M^\circ_{F'}G'/A_{*+1}C_{F'}G'
\stackrel{s^*}\to
A_*C_{F'}G' \to A_*(C_FG\times_FF')
\end{align*}
(the last map is pushforward via the closed immersion
$C_{F'}G'\to C_FG\times_FF'$.
In case $f\colon F\to G$ is a regular local immersion of
codimension $d$,
the cone $C_FG$ is a bundle, and we postcompose with
the inverse to the pullback map
$A_{*-d}F'\eqto A_*(C_FG\times_FF')$
to obtain, finally,
the {\em refined Gysin homomorphism} $f^!\colon A_*G'\to A_{*-d}F'$.

The refined Gysin map $f^!\colon A_*G'\to A_{*-d}F'$, on the level
of actual cycles, amounts
to nothing more than the usual refined Gysin map.
If $T\to G'$ is a projective morphism and $E\to T$ is a vector bundle,
then applying $f^!$ to the element of $A_*G'$ given by
$[V]\in A^\circ_*E$ yields the cycle class in $A_*F'$ given by
$f^!([V])\in A^\circ_*(E\times_FG)$.

By purely formal arguments we see that the refined Gysin map to
a regular local immersion is compatible with flat pullback and
with projective pushforward.
Functoriality of the refined Gysin map follows
exactly as in \cite[\S6.5]{f}
(the claim for general $F\to G\to H$ reduces to the case
when $H$ is a vector bundle over $G$ and $G\to H$ is the zero section,
in which case the result follows from a local computation showing
$N_FH\simeq N_FG\oplus H|_F$).
The argument for commutativity of refined Gysin homomorphisms corresponding to
a pair of regular local immersions of stacks reduces by formal
manipulations to the statement of \cite[Theorem 6.4]{f},
and this we know by \cite[Proposition 4]{kr}.

\subsection{Exterior product}
\label{exterior}

Let $X$ and $Y$ be stacks.
Since the product of projective morphisms is again projective,
a pair consisting of a cycle
$(f,\alpha)$, with $f\colon S\to X$ and $\alpha\in \hat A_*S$,
and a cycle
$(g,\beta)$, with $g\colon T\to Y$ and $\beta\in \hat A_*T$,
determines a cycle $(f\times g, \alpha\times \beta)\in A_*(X \times Y)$.

\begin{pr}
The map, sending $((f,\alpha), (g,\beta))$ to
$(f\times g, \alpha\times \beta)$,
determines a morphism $A_*X\otimes A_*Y\to A_*(X\times Y)$.
\end{pr}

The proof involves only routine checking of details.

\subsection{Intersections on Deligne-Mumford stacks}

Intersection theory on stacks was motivated by a desire to
establish foundations for enumerative calculus on
moduli spaces.
A particularly enlightening early investigation in this
direction is \cite{m}, where foundations are laid for
constructing intersection rings with rational coefficients
on certain moduli spaces of curves.
The intersection ring of the compactified moduli space $\overline M_2$
of curves of genus 2 is described in detail,
with the fractional coefficients that appear attributed to
automorphisms of the curves.
Later, more general, approaches to intersection theory
on Deligne-Mumford stacks \cite{g,v}
also produce intersection operations and intersection products
on the Chow groups with rational coefficients.

The functor $A_*$ allows the construction of an integer-valued
intersection product.
If $X$ is a smooth Deligne-Mumford stack, then the diagonal
$X\to X\times X$ is a regular local immersion,
so the contents of sections \ref{fmp} and \ref{exterior} yield
an intersection product on $A_*X$.
This intersection product satisfies the usual properties,
cf.\ \cite[\S8.3]{f}.
In particular, it agrees with the intersection product of
\cite{eg} in case $X$ is a global quotient.

Even though we have an integer-valued intersection ring,
we still need to tensor with $\Q$
if we wish to do enumerative geometry.
If $X$ is a complete Deligne-Mumford stack over a field $k$, then
the map $X\to \Spec k$ is proper
(cf.\ \cite{dm} or \cite{v} for a definition)
but non-representable, so we do not get a pushforward in the $A_*$ theory.
What we have is a cycle map back to the na\"\i{}ve Chow groups
$$cyc\colon
A_*X\to A_*X\otimes\Q\stackrel{\raise-1.2pt\hbox{$\scriptstyle\sim$}}\to
A^\circ_*X\otimes\Q,$$
which introduces denominators from the fact that
the Gysin map for vector bundles is defined only rationally on $A^\circ_*$.
Then, we have a pushforward on na\"\i{}ve Chow groups
$$\textstyle \int_X\colon A^\circ_0X\otimes\Q\to A^\circ_0\Spec k\otimes\Q\simeq\Q,$$
cf.\ \cite{v},
and this pushforward may introduce even more denominators.
As an example, the compactified moduli space
$\overline M_{1,1}$ of elliptic curves
over the complex numbers has two special points with
stabilizer groups cyclic of orders 4 and 6, respectively,
and generic point with stabilizer group $\Z/2$.
For all $\alpha\in A_0 \overline M_{1,1}$ we have
$2\,cyc(\alpha)\in \im(A^\circ_0 \overline M_{1,1}\to
A_0 \overline M_{1,1})$.
However, if we let
$\pi\colon U\to \overline M_{1,1}$ be the universal curve and
let $E=\pi_* \omega_{U/\overline M_{1,1}}$ be the Hodge bundle,
then we find $\int_{\overline M_{1,1}} cyc(c_1(E))=1/24$.

\subsection{Boundedness by dimension}

Since projective pushforward lowers codimension,
there is the potential, a priori, that there can be
nontrivial cycle classes in $A_kX$ for $k$ greater
than the dimension of $X$.
This turns out not to occur; the justification
uses some facts about projective morphisms
plus the splitting principle.
We first prove a preliminary lemma, and then prove
the vanishing of $A_kX$ for $k>\dim X$.
This section, as well as our deduction of
Corollary \ref{xtendproj} from Proposition \ref{xtendcoh},
make clear why we do
not develop the pushforward for general proper morphisms, but only
for projective morphisms.

\begin{lm}
\label{vanishp}
Let $Y$ be a stack of dimension $d$,
let $E\to Y$ be a vector bundle of rank $r$ with
projectivization $p\colon P(E)\to Y$,
and let $U_1$ and $U_2$ be vector bundles on $Y$
of ranks $e_1$ and $e_2$, respectively.
Suppose $\gamma\in A^\circ_k((p^*U_1)(m)\oplus p^*U_2)$
with $k>d+e_1+e+2$
(here $(p^*U_1)(m)$ denotes
$(p^*U_1)\otimes \cO_E(1)^{\otimes m}$).
Then $(p,\gamma)=0$ in $A_*Y$.
\end{lm}

\begin{proof}
We induct on $r$.
The case $r=1$ is trivial.
For the inductive step,
we consider
$q\colon P:=P(E^\vee)\to Y$,
with $q^*E\to \cO_{E^\vee}(1)\to 0$ on $P$.
If we let
$L=\cO_{E^\vee}(-1)$ then we have
$$0\to K\to q^*E\otimes L\to \cO_P\to 0$$
on $P$,
for some $K$ of rank $(r-1)$.
Let $E'=p^*E\otimes L$.
Then $\cO_{E'}(1)$ has a global section which vanishes
precisely on $P(K)\subset P(E')$.
Moreover we may identify $P(E')$ with $P(q^*E)$, i.e.,
$P(E')$ fits into a fiber diagram
$$\xymatrix{
P(E') \ar[r]^{p'} \ar[d]_{q'} & P \ar[d]^q \\
P(E) \ar[r]^p & Y
}$$
and moreover if we make this identification then we find
$q^{\prime{*}}((p^*U_1)(m))\simeq (p^{\prime{*}}U'_1)\otimes
\cO_{E'}(1)^{\otimes m}$,
where $U'_1=q^*U\otimes L^{\otimes m}$.
If we set $U'_2=q^*U_2$ then we have
$q^{\prime{*}}\simeq p^{\prime{*}}U'_2$,
so via these identifications we have
$$q^*(p,\gamma)=(p',\delta)$$
with $\delta\in A^\circ_{k+r-1}((p^{\prime{*}}U'_1)(m)\oplus
p^{\prime{*}}U'_2)$.

Since $\cO_{E'}(1)$ has a section nonvanishing on
$P(E')\setminus P(K)$,
we may find
$\delta'\in A^\circ_{k+r-1}(p^{\prime{*}}U'_1\oplus p^{\prime{*}}U'_2)$
such that $\delta$ and $\delta'$ have the same image under
restriction to
$A^\circ_{k+r-1}(u^*p^{\prime{*}}U'_1\oplus u^*p^{\prime{*}}U'_2)$,
where
$u\colon P(E')\setminus P(K)\to P(E')$ denotes inclusion.
If we consider the projections
$$\xymatrix{
& {(p^{\prime{*}}U'_1)(m)\oplus (p^{\prime{*}}U'_2)} \\
{(p^{\prime{*}}U'_1)(m)\oplus (p^{\prime{*}}U'_1) \oplus (p^{\prime{*}}U'_2)}
\ar[ur]^{pr_{13}}\ar[dr]_{pr_{23}} \\
& {(p^{\prime{*}}U'_1)\oplus (p^{\prime{*}}U'_2)}
}$$
we have, by Remark \ref{twovbsurj} plus the standard
excision sequence for na\"\i{}ve Chow groups,
$$pr^*_{13}\delta = pr^*_{23}\delta'+i'_*\varepsilon$$
for some
$\varepsilon\in A^\circ_{k+r-1+e_1}((i^*p^{\prime{*}}U'_1)(m)\oplus
(i^*p^{\prime{*}}U'_1)\oplus (i^*p^{\prime{*}}U'_2))$,
where $i\colon P(E)\to P(E')$ denotes inclusion
and $i'$ denotes the pullback of $i$.
We find $(p',\delta')=0$ for dimension reasons,
and $(p'\smallcirc i,\varepsilon)=0$ by the induction hypothesis,
and hence $(p',\delta)=0$ in $A_*P$.
Since $q^*$ is injective, we have
$(p,\gamma)=0$ in $A_*Y$ as desired.
\end{proof}

\begin{pr}
Let $Y$ be a stack.
We have $A_k(Y)=0$ for all $k>\dim Y$.
\end{pr}

\begin{proof}
It suffices to show that for any
$\alpha\in A_k(Y)$ that there exists a nonempty open
substack $U$ of $Y$ such that,
with inclusion map $\nu\colon U\to Y$,
we have $\nu^*\alpha=0$ in $A_k(U)$.
It suffices to consider $\alpha=(f,\gamma)$ with $f\colon T\to Y$
projective and $\gamma\in A^\circ_*F$ with $F$ a vector bundle over $T$.
Shrinking $Y$, we may assume $f$ factors (up to 2-isomorphism) as
a closed immersion followed by a projection of the form
$P(E)\to Y$ where $E$ is a vector bundle on $Y$.
There is a bundle $\cO(1)$ on $P(E)$, hence also on $T$,
and if we let $\cF$ denote the sheaf of sections of $F$,
then the natural map $f^*f_*\cF(m)\to \cF(m)$ is surjective
for suitable $m$.
Shrinking $Y$ further, we may suppose $f_*\cF(m)$ is locally free,
so we have a surjection of vector bundles
$(f^*V)(-m)\to F$,
and so we are reduced to the case $\alpha=(p,\gamma)$ with
$p\colon P(E)\to Y$ the projection map and
$\gamma\in A^\circ_*((p^*V)(-m))$ for $V$ a vector bundle on $Y$.
But now $(p,\gamma=0)$ in $A_*Y$ by Lemma \ref{vanishp},
so we are done.
\end{proof}

\subsection{Stratifications by quotient stacks}

As promised in the introduction,
we shall eventually show that whenever $\pi\colon B\to Y$
is a vector bundle stack, such that the base $Y$ admits a stratification
by global quotient stacks,
then $\pi^*\colon A_*Y\to A_*B$ is an isomorphism.
To prove this requires the localization machinery of section \ref{xxa}.
At this point, we content ourselves with some elementary observations,
first describing classes of stacks which admit such stratifications
and then showing (Proposition \ref{schemeapprox} (ii))
how to obtain the desired homotopy property when
$Y$ is a suitable global quotient
(this includes many cases of interest,
e.g., certain moduli spaces).

\begin{lm}
\label{quotientlm}
Let $X$ be an algebraic space, and let
$G$ be a linear algebraic group acting on $X$.
Then $[X/G]^{\rm red}$ contains a nonempty open substack
isomorphic to $[V/GL(n)]$ for some
quasiprojective scheme $V$ with linear action of
$GL(n)$.
\end{lm}

\begin{proof}
Replacing $X$ by $X\times GL(n)/G$ we may assume $G=GL(n)$ for some $n$;
then $[X/G]^{\rm red}=[X^{\rm red}/G]$, so we may as well assume $X$
reduced.
Choosing an irreducible affine open $U\subset X$, we may replace $X$
by the image of the action map $U\times G\to X$,
and now there must exist a finite subset $S=\{g_i\}$ of closed points of $G$
such that $U\times S\to X$ is surjective and
such that the residue field of each $g_i$ is separable over the base field.
So, for some finite separable field extension $k\to k'$ of the base field,
there is a cover of $X_{k'}:=X\times_{\Spec k}\Spec k'$ by affine schemes,
and hence $X_{k'}$ is a scheme.
By \cite{su}
there exists a quasiprojective $G$-stable dense open
subscheme
$V\subset X_{k'}$.
If $Z$ denotes the complement $X_{k'}\setminus V$ then the
image $Y$ of $Z$ under the finite map $X_{k'}\to X$ is a
proper closed subscheme of $X$, and now $X\setminus Y$ is
$G$-stable and is a quasiprojective scheme
(since $(X\setminus Y)_{k'}$ is a quasiprojective scheme).
Finally, now,
$[(X\setminus Y)^{\rm reg}/G]$ is a nonempty open substack of our
original stack which is the quotient by $G=GL(n)$ of a
quasiprojective scheme that is regular, and in particular normal,
so the action of the group is linearizable.
\end{proof}

\begin{pr}
\label{stratprop}
Let $Y$ be a stack.
The following are equivalent.
\begin{romanlist}
\item There exists a stratification of $Y^{\rm red}$
by locally closed substacks $U_i$
such that each $U_i$ is isomorphic to a stack of the form
$[X_i/G_i]$, where for each $i$, $X_i$ is an algebraic space,
and $G_i$ is a linear algebraic group acting on $T_i$;
\item There exists a stratification of $Y^{\rm red}$
by locally closed substacks $U_i$
such that each $U_i$ is isomorphic to a stack of the form
$[T_i/G_i]$, where for each $i$, $T_i$ is a quasiprojective scheme
and $G_i$ is a smooth connected linear algebraic group acting linearly
on $T_i$.
\end{romanlist}
\end{pr}

\begin{proof}
Immediate from Lemma \ref{quotientlm}.
\end{proof}

\begin{defn}
A stack $Y$ is said to admit a {\em stratification by global quotients}
if the conditions of Proposition \ref{stratprop} are satisfied for $Y$.
\end{defn}

\begin{conv}
The words {\em global quotient} or {\em quotient stack},
without any additional qualifyers,
refer from now on to quotients of an algebraic space by an arbitrary linear
algebraic group, as in (i) above.
\end{conv}

\begin{pr}
\label{gq}
\begin{romanilist}
\firstitem Let $X$ and $Y$ be stacks which admit stratifications
by global quotients.
Then $X\times Y$ admits a stratification by global quotients.
\item Let $Y$ be a stack which admits a stratification by global quotients,
and let $f\colon X\to Y$ be a representable morphism.
Then $X$ admits a stratification by global quotients.
\item Every Deligne-Mumford stack admits a stratification by
global quotients.
\end{romanilist}
\end{pr}

\begin{proof}
Since the product of global quotient stacks
is again a global quotient stack, (i) is clear.
Claim (ii) follows from the fact that for any algebraic group $G$,
a representable morphism $U\to BG$ leads to an action of $G$
on $X:=U\times_{BG}\Spec k$ such that $U\simeq [X/G]$.
For (iii), if $f\colon U\to F$ is an \'etale presentation of a
Deligne-Mumford stack $F$, then the restriction of $f$ over some
nonempty open substack $G$ of $F$ is finite \'etale of some degree $n$,
and $G$ is isomorphic to the quotient of the complement of
all the diagonal components of $U\times_GU\times_UG\cdots\times_GU$
($n$ copies) by the symmetric group $S^n$
(cf.\ \cite[(10.2)]{l}).
\end{proof}

\begin{pr}
Let $X$ be a stack.
The following are equivalent.
\begin{romanlist}
\item For every integer $N$,
there exist a vector bundle $E\to X$ and a representable open substack
$U$ of $E$ such that $E\setminus U$ has codimension $\ge N$ in $E$.
\item There exist a vector bundle $E\to X$
and a locally closed immersion $T\to E$,
with $T$ representable and $T\to E$ surjective.
\item There exists an algebraic space $P$ with action of $GL(n)$
for some $n$, such that
$X\simeq [P/GL(n)]$.
\item $X$ is a global quotient stack.
\end{romanlist}
\end{pr}

\begin{proof}
Conditions (iii) and (iv) are equivalent and, by the construction of \cite{eg},
imply (i).
Clearly (i) implies (ii).
Suppose (ii) holds, and let $P$ be the principal bundle associated
to the vector bundle $\pi\colon E\to X$, so that
$X\simeq [P/GL(n)]$,
where $n=\rk E$.
We claim $P$ must be representable.

We may assume $k$ is algebraically closed.
We may also assume $T$ is disjoint from the zero section
$s(X)$ (replace $T$ by $[T\cup \pi^{-1}(T\cap s(X))]\setminus s(X)$).
Choosing, say, the first basis element of a framing yields a
representable, faithfully flat morphism $P\to E\setminus s(X)$.
The pre-image of $T$ is a representable,
locally closed substack $S$ of $P$,
such that the translates of $S$ by elements of $GL(n)$
cover the $k$-valued points of $P$.
Hence $P$ is representable.
\end{proof}

When the base field has positive characteristic, there exist stacks $Y$
which have finite stabilizer at every point but
are not Deligne-Mumford.

\begin{pr}
\label{isaquotient}
Let $Y$ be a stack.
Then the following are equivalent.
\begin{romanlist}
\item The diagonal map $Y\to Y\times Y$ is quasi-finite.
\item The stabilizer $Y\times_{Y\times Y}Y\to Y$ is quasi-finite.
\end{romanlist}
Moreover, if $Y$ has quasi-finite diagonal then $Y$ admits a stratification
by global quotients.
\end{pr}

\begin{proof}
Clearly (i) implies (ii).
For the converse, it suffices to check that whenever $\Omega$ is an
algebraically closed field containing the base field $k$,
then for any $x,y\in Y(\Omega)$
the set $\Isom_Y(x,y)(\Omega)$ is finite.
But if $\Isom_Y(x,y)(\Omega)$ is nonempty,
then for any $t\in \Isom_Y(x,y)(\Omega)$,
postmultiplication with the inverse to $t$
gives an injective map $\Isom_Y(x,y)\to \Isom_Y(x,x)$,
so (ii) implies (i) (cf.\ \cite{dm} or \cite{v} for notation).
Now assume $Y$ has quasi-finite diagonal.
We replace $Y$ by $T^{\rm red}$
and let $f\colon U\to Y$ be a smooth presentation.
By \cite[(5.7)]{l},
for a suitable closed subscheme $V$ of $U$,
the map $V\to Y$ is dominant,
and the restriction over a dense open substack of $Y$
is finite and flat.
Let us shrink $Y$;
if $E\to Y$ denotes the vector bundle whose sheaf of sections is
$f_*\cO_V$,
then we have a closed immersion $V\to E$,
so by Proposition \ref{isaquotient}, $Y$ is a global quotient.
\end{proof}

\begin{pr}
\label{schemeapprox}
Let $X$ be a quasiprojective scheme, let $G$ be a connected
smooth linear algebraic group acting linearly on $X$,
and let $Y=[X/G]$.
Suppose $\pi\colon B\to Y$ is a vector bundle stack.
Then
\begin{romanlist}
\item $B$ has vector bundles with total spaces represented by
schemes off of loci of arbitrarily high codimension; and
\item the map $\pi^*\colon A_*Y\to A_*B$ is an isomorphism.
\end{romanlist}
\end{pr}

\begin{proof}
Since $G$ acts linearly on the quasiprojective scheme $X$,
the quotient stack $Y$ has the property
that every coherent sheaf admits a surjective map from a
locally free sheaf.
So, by \cite[Proposition 1.4.15]{d},
$B$ is isomorphic to a globally presented vector bundle stack
$[E/F]$ for some morphism $F\to E$ of vector bundles on $Y$.
To establish (i) and (ii), it suffices to consider the case
when $E$ is zero, i.e., $B=BF$.
Now we consider the vector bundle $R\to B$ given by the bundle
$F\oplus 1$ on $Y$ with $F$-action given by
$f\colon (x,t)\mapsto (x+tf,t)$.
There is a projection map $r\colon R\to\A^1$,
and the fiber $r^{-1}(\A^1\setminus\{0\})$ is isomorphic
to $Y\times(\A^1\setminus\{0\})$.
Hence the composite map
$R^{\oplus n}\to B\to Y$,
restricted to $r^{\times n}(\A^n\setminus\{0\})$,
is representable,
and so for a suitable vector bundle $E\to B$
(we take $E$ to be the pullback from $BG$ of a suitable
representation bundle of $G$,
as in \cite{eg}),
the total space of $R^{\oplus n}\oplus E$ possesses an
open substack that is representable by a scheme
and has complement of codimension $\ge n$.
This establishes (i).
Statement (ii) is a consequence of Corollary \ref{pistarsplitmono},
since the structure map $\varphi\colon Y\to B$ is an affine bundle
and we have $\pi\smallcirc \varphi = 1_Y$.
\end{proof}

\section{Extended excision axiom}
\label{xxa}

\subsection{A first higher Chow theory}

We need some sort of a first higher Chow group in order to
be able to extend
the excision sequence one place to the left.
We take as motivation the long exact localization sequence coming from
the Gersten complex of schemes.

\begin{nota}
In this section, $A_j(X;1)$ denotes, for a scheme $X$,
the kernel of $\partial\colon W_jX\to Z_jX$, modulo the
subgroup generated by tame symbols of elements of
$K_2(k(Y))$ for all $(j+2)$-dimensional integral closed subschemes $Y$ of $X$.
\end{nota}

\begin{rem}
On a separated scheme $X$,
we recognize $A_j(X;1)$ as an $E^2$ term of
the Quillen spectral sequence of $K$-theory.
Formal properties, cf.\ \cite{gadv}, imply that if
$X$ is a scheme (separated or not), and if
$\pi\colon E\to X$ is a flat morphism whose pointwise fibers are
$r$-dimensional affine spaces, then the induced map
$A_j(X;1)\to A_{j+r}(E;1)$ is an isomorphism.
\end{rem}

The failure of descent for
$K$-theory in the smooth topology
means we cannot apply the machinery of the Quillen spectral sequence
directly to stacks.
We must resort to a cycle-based complex, which we can show
to be quasi-isomorphic (in the needed range)
to the Gersten complex on a scheme base.

The starting point for cycle-based complexes is the
theory of Bloch's higher Chow groups \cite{badv,bjag}.
Let us recall that the $n^{\rm th}$ term in the Bloch complex for a scheme $X$
is the free abelian group of cycles
on $X\times \triangle^n$ meeting boundary cycles properly;
$\triangle^n$ denotes the algebraic $n$-simplex ($\simeq \A^n$)
and the boundary cycles
are copies of $X\times \triangle^m$ for $m<n$.
In particular, when $n=1$ the boundary consists of just two
zero-simplices (points), and the condition of proper intersection
says nothing other than that no component of a cycle be contained
in either of the boundary components.
Thus the rightmost terms in the Bloch complex are
\begin{equation}
\label{blochright}
\cdots\to Z_{j+1}X\times(\A^1\setminus\{\rm 2\ points\})\to
Z_jX\to 0,
\end{equation}
where the final boundary map is the difference of two cycle-level
specialization maps.

We identify
$\A^1\setminus\{\rm 2\ points\}$
with $R:=\bP^1\setminus\{0,-1,\infty\}$,
and we denote by $\pi$ the projection $X\times R\to X$.
For $t\in\bP^1$, we define
$\partial_t\colon Z_{j+1}X\times\bP^1\to Z_jX$
to be the cycle-level pullback via $X\simeq X\times\{t\}\to X\times\bP^1$.
Since $\partial_t$ kills any cycle supported in a fiber of $\pi$,
there is an induced map
$\partial_t\colon Z_{j+1}X\times R\to Z_jX$.
We set $\partial=\partial_0-\partial_\infty$.
Then $\partial\colon Z_*X\times R\to Z_*X$ is the
rightmost map in (\ref{blochright}).

We remark that there is a variant of (\ref{blochright}),
obtained by moding out by degenerate cycles
(so, e.g., in term 1 we mod out by $\pi^*Z_jX$).
The so-called normalized complex which results is quasi-isomorphic to
the original complex.

Let $X$ be a stack.
We introduce a cycle complex which has the same groups in positions
0 and 1 as the normalized Bloch complex, and we put the group
$Z_*(X\times T^2)$ in position 2 to yield
\begin{equation}
\label{adhoccx}
Z_{j+2}(X\times T^2)\stackrel{\partial}\to
Z_{j+1}(X\times R) / \pi^*Z_jX\stackrel{\partial}\to
Z_jX\to 0.
\end{equation}
Here $T^2$ denotes the two-dimensional torus $(\A^1\setminus\{0\})^2$.
The complex will not go beyond position 2.

The map
$\partial\colon Z_{j+1}(X\times R) / \pi^*Z_jX\to Z_jX$
is the map $\partial_0 - \partial_\infty$ of the Bloch complex.
We give a definition of the boundary map
$\partial\colon Z_{j+2}(X\times T^2)\to Z_{j+1}(X\times R) / \pi^jZ_*X$,
so that (\ref{adhoccx}) is a complex.
We start by fixing an orientation convention on toric compactifications
of $T^2$.
Let $Y$ be a nonsingular two-dimensional complete toric variety,
with corresponding fan $\Delta$.
Let $(u,v)\in N=\Z^2$ be a generator of a ray $\rho\in\Delta$.
Let $(u',v')$ be the generator of the ray immediately preceding
$\rho$ via the counterclockwise ordering of rays, and let
$(u'',v'')$ generate the ray immediately following $\rho$.
If $x$ and $y$ denote coordinates on $T^2\subset Y$,
then the maximal cone of $\Delta$ preceding $\rho$ corresponds
to the toric affine chart
$\Spec k[x^vy^{-u},x^{-v'}y^{u'}]$, and the maximal
cone of $\Delta$ following $\rho$ corresponds
to the toric affine chart
$\Spec k[x^{v''}y^{-u''},x^{-v}y^u]$.
Corresponding to $\rho$ is a toric divisor
$D_\rho$ of $Y$.
We define our orientation convention to be the
identification of $D_\rho$ with $\bP^1$ via
$\Spec k[x^vy^{-u}]\to \A^1=\bP^1\setminus\{\infty\}$
and $\Spec k[x^{-v}y^u]\to \bP^1\setminus\{0\}$.
There are two possible orientation conventions,
but the point of $D_\rho$
corresponding to $\{-1\}$ in $\bP^1$ is independent of this choice.
Thus there is a natural subset
$D_\rho^\circ\subset D_\rho$,
defined as the complement in $D_\rho$ of the torus-fixed points and
the point corresponding to $\{-1\}$.
Our convention specifies an isomorphism $R\to D_\rho^\circ$.

Suppose $V$ is an integral closed substack of $X\times\bP^2$ such that
$V$ meets $X\times T^2$ nontrivially,
By induction on excess of intersection (cf.\ \cite[\S2.4]{f}),
we find that there is a finite sequence of blowups
at torus-fixed points $Y\to \bP^2$ such that the
proper transform $\widetilde V$ of $V$ meets the toric divisors
of $X\times Y$ properly.
For each $\rho$ in the fan defining $Y$ we can pull back
$[\widetilde V]$ via the composite
\begin{equation}
\label{mapfromr}
X\times R \eqto X\times D_\rho^\circ\to X\times D_\rho\to X\times Y
\end{equation}
to obtain a cycle $\partial_\rho([V])\in Z_*X\times R$.
By dimension reasoning,
the sum
$\sum_\rho \partial_\rho([V])$ in $Z_*(X\times R)/\pi^* Z_*X$ is independent
of the choice of $Y$.
Because of proper intersection, we have
$\partial(\sum_\rho \partial_\rho([V]))=0$.
Hence, if we define $\partial=\sum_\rho \partial_\rho$, then
(\ref{adhoccx}) is a complex.

\begin{defn}
Let $X$ be a stack.
We denote by $\underline A^\circ_jX$ the homology group in
the first position of the complex (\ref{adhoccx}).
\end{defn}

\begin{rem}
The boundary maps clearly respect proper pushforward and flat pullback,
making the association
$X\mapsto\underline A^\circ_*X$
functorial for proper pushforward and flat pullback.
\end{rem}

Let us examine the case $X=\Spec k$.
Let $C\subset T^2$ be an integral subscheme of dimension 1.
Then $C$ is given as the zero locus of a single function
$f\in k[x^{\pm 1},y^{\pm 1}]$.
It is easy to describe a toric variety $Y$ suitable for computing
$\partial([C])$.
We let $N=\Z^2$, and $M=\Hom(N,\Z)$.
If $f=\sum a_{\mu\nu} x^\mu y^\nu$, we let
$\Gamma$ be the Newton polygon of $f$,
that is,
the convex hull in $M\otimes\R$ of the set of points $(\mu,\nu)$
for which $a_{\mu\nu}\ne 0$.
There exists a finite collection of half-planes
$\{H^+_{\rho_i}\}$ which cut out $\Gamma$,
where for $\rho\in N$ we define
\begin{align*}
\lambda_\rho &= \min\{\,\langle \zeta,\rho\rangle\mid
\zeta\in\Gamma\,\}, \\
H_{\rho}&=
\{\,\zeta\in M\otimes\R\mid \langle \zeta,\rho\rangle=\lambda_\rho\,\}, \\
H^+_{\rho}&=\{\,\zeta\in M\otimes\R\mid \langle \zeta,\rho\rangle \ge
\lambda_\rho\,\},
\end{align*}
where $\langle\,,\,\rangle$ denotes the pairing of $M$ and $N$.
If $p$ and $q$ lie in the set $\VVert(\Gamma)$
of vertices of $\Gamma$, let us say that
$q$ follows $p$ if
the line segment joining $p$ and $q$ is an edge of $\Gamma$,
and if for all $r\in\intrior(\Gamma)$,
the signed angle from ray $pq$ to ray $pr$ is positive.

\begin{pr}
\label{whichfan}
Suppose $\Gamma=\bigcap_i H^+_{\rho_i}$.
If $\Delta$ is a complete nonsingular fan in $N$
which contains all the rays generated by the $\rho_i$,
then the closure $\widetilde C$ of $C$ in the toric variety
corresponding to $\Delta$ meets the toric divisors properly.
Moreover, if $N\colon Z_0 R\to k^*$ is defined by
sending a zero-cycle $[Z]$, with $Z\subset R$ integral,
to the image under the norm map $k(Z)^*\to k^*$ of the function
$(-t)$, where $t$ is the restriction to $R$ of the natural coordinate
on $\A^1$, then we have
$$N(\partial_\rho([C]))=
\begin{cases}
a_p/a_q & \text{if
$H_{\rho}\cap\VVert(\Gamma)=\{p,q\}$ such that $q$ follows $p$,} \\
1 & \text{otherwise.}
\end{cases}
$$
\end{pr}

\begin{proof}
Proper intersection is equivalent to saying that
$\widetilde C$ does not contain any of the torus fixpoints.
Let $\rho_1$, $\ldots$, $\rho_m$ be generators of the rays
of $\Delta$, arranged in counterclockwise order.
The hypotheses guarantee that $H_{\rho_i}$ and $H_{\rho_{i+1}}$
intersect at a vertex of $\Gamma$, for each $i$.
Thus the equation for $\tilde C$ in a neighborhood of
the $i^{\rm th}$ torus fixpoint has a nontrivial constant term.
Moreover, $\widetilde C\cap D_{\rho_i}\ne\emptyset$ if and only if
$H_{\rho_i}$ contains two vertices of $\Gamma$,
say $p=H_{\rho_{i-1}}\cap H_{\rho_i}$ and
$q=H_{\rho_i}\cap H_{\rho_{i+1}}$.
For such $i$,
the equation defining the scheme $C\cap D_i$ is a polynomial
with leading term $a_q$ and constant term $a_p$,
and hence the formula.
\end{proof}

\begin{cor}
\label{isospeck}
The composite $Z_1T^2\to Z_0R\to k^*$ is the zero map, and
the induced map $Z_0R/\partial(Z_1T^2)\to k^*$ is an isomorphism.
\end{cor}

\begin{proof}
The first claim is an immediate consequence of
Proposition \ref{whichfan}.
For the second claim, the map is clearly surjective,
so we need only verify that a general zero-cycle on $R$ is
equal, modulo $\partial$ of an element of $Z_1T^2$,
to a cycle of the form $[\{r\}]$ for some $r\in k^*$.
A general effective zero-cycle is of the form $[Z]$
where $Z$ is the zero locus of some polynomial
$t^n+a_1t^{n-1}+\cdots+a_n$.
If $C$ is the zero locus of
$x^n+a_1x^{n-1}+\cdots+a_n+y$,
then $\partial([C])=[Z]+[\{r^{-1}\}]$, where $r=-a_n$.
In particular, as well, $[\{r\}]+[\{r^{-1}\}]=0$.
Thus the class of a general zero-cycle has a representative
of the form $[\{r\}]$ for some $r\in k^*$.
\end{proof}

\begin{pr}
\label{nqiso}
There is a natural isomorphism of functors on schemes
$\underline A^\circ_*({-})\to A_*({-};1)$.
\end{pr}

\begin{proof}
Let $X$ be a scheme, and let $j$ be an integer.
We propose maps $N$ yielding a morphism of complexes
\begin{equation}
\label{normmap}
\begin{split}
\xymatrix{
Z_{j+2}(X\times T^2) \ar[r]^(.4)\partial \ar[d]_{N} &
Z_{j+1}(X\times R)/\pi^*Z_jX \ar[r]^(.7)\partial \ar[d]_{N} &
Z_jX \ar[d]^{{\rm id}} \\
{\coprod_{x\in X_{j+2}} K_2(k(x))} \ar[r] & W_jX \ar[r] & Z_jX
}
\end{split}
\end{equation}
from the cycle complex (\ref{adhoccx}) to the Gersten complex of $X$
(we denote by $X_j$ the set of points $x\in X$ having
$\dim \overline{\{x\}}=j$).
The middle vertical map is the norm map, preceded by the
involution $t\mapsto -t$ on $R$ (as in Proposition \ref{whichfan}).
For the map on the left,
let $V$ be an integral closed subscheme of $X\times\bP^2$
of dimension $(j+2)$,
and consider the image
$p(V)$ under projection $p\colon X\times \bP^2\to X$.
There are 3 cases:
\begin{romanlist}
\item $\dim p(V)=j$, i.e., $V=U\times\bP^2$ for some $U\subset X$,
and clearly $[V]$ is killed by $\partial$; we define $N([V])=0$.
\item $\dim p(V)=j+1$, so $N\smallcirc \partial([V])$ consists just of a
rational function on $p(V)$.  To compute, we may replace $X$ by the
generic point of $p(V)$, and by Corollary \ref{isospeck},
$N(\partial([V]))=0$; we define $N([V])=0$.
\item $\dim p(V)=j+2$, so $k(p(V))\to k(V)$ is a finite field extension,
and there is a norm map $K_2(k(V))\to K_2(k(p(V)))$.
We send $[V]$ to the image under the norm map of the
symbol $\{-x,-y\}\in K_2(k(V))$,
where $x$ and $y$ denote the coordinate functions on $T^2$.
\end{romanlist}

We have defined a map
$N\colon Z_{j+2}(X\times\bP^2)\to \coprod_{x\in X_{j+2}} K_2(k(x))$.
Since this map vanishes on cycles supported in the complement of $X\times T^2$,
there is an induced map $N$ on $Z_{j+2}(X\times T^2)$,
as indicated in (\ref{normmap}).
Now, commutativity of (\ref{normmap}) follows from the fact that the
norm map of $K$-theory commutes with the tame symbol.

We claim, now, that $N$ induces
isomorphisms on the zeroth and first homology groups.
The map on zeroth homology groups is clearly the identity map on $A^\circ_*X$.
Let us consider the induced map $N_1$ on first homology groups.
The vertical maps $N$ in (\ref{normmap}) are surjective
(given any rational function $r$ on an integral
closed subscheme of $X$, the graph of $(-r)$ specifies an element of
$Z_*(X\times R)$ whose image under $N$ is the specified rational function;
a similar argument with pairs of functions applies for the map on the left),
so $N_1$ is surjective, and to show $N_1$ is injective it suffices to
show that if $\alpha\in Z_*(X\times R)/\pi^*Z_*X$ satisfies
$N(\alpha)=0$ in $W_*X$, then $\alpha$ lies in the image of $\partial$.
For this we are easily reduced to the case $X=\Spec k$,
and by Corollary \ref{isospeck},
$N_1\colon \underline A^\circ_*(\Spec k)\to A_*(\Spec k;1)$
is an isomorphism.
\end{proof}

We wish to repeat the construction of section \ref{homologyfunctor},
starting with the functor $\underline A^\circ_*$ rather than $A^\circ_*$.
Thus we need to study what happens when we have two surjections of
vector bundles $E\to F$ on a stack $X$.

\begin{lm}
Let $S=T^2\setminus\{\,(x,y)\mid x+y+1=0\,\}$.
The map $\partial\colon Z_*X\times T^2\to Z_*(X\times R)/\pi^*Z_*X$
factors through $Z_*(X\times S)$.
The map $\partial$ also factors through $Z_*(X\times R\times R)$.
\end{lm}

\begin{proof}
Suppose $\alpha\in Z_*X\times T^2$ has support in the locus specified by
$x+y+1=0$.
Then the closure of $\alpha$ in $Z_*X\times\bP^2$ meets the boundary
divisors properly, and the pullback of $\alpha$ under
(\ref{mapfromr}) is zero for each of the 3 boundary cycles of $\bP^2$.
A similar argument applies if $\alpha$ has support in the locus
given by $(x+1)(y+1)=0$, using
$\bP^1\times\bP^1$ in place of $\bP^2$.
\end{proof}

There is a map $S\to R$, given by
$(x,y)\mapsto x+y$.
If $\sigma$ denotes the induced map
$X\times S\to X\times R$, then for any
$\alpha\in Z_*(X\times R)$,
$\sigma^*\alpha$ satisfies $\partial(\sigma^*\alpha)=\alpha+\tau^*\alpha$,
where $\tau$ is the map
$X\times R\to X\times R$ induced by the involution $t\mapsto t^{-1}$
on $R$.

\begin{pr}
\label{twovbua}
Let $X$ be a stack, let $E$ and $F$ be vector bundles on $X$,
and let $\varphi$ and $\psi$ be surjections of vector bundles
from $E$ to $F$.
Then the maps $\underline A^\circ_*F\to \underline A^\circ_*E$
induced by $\varphi$ and $\psi$ are the same.
\end{pr}

\begin{proof}
Let $r$ and $s$ denote free parameters;
then $r\varphi+s\psi$ is a vector bundle surjection
from $E\times \A^2$ to $F$.
Consider the restriction of
$r\varphi+s\psi$ to
$E\times(\A^2\setminus P)$,
where $P=\{\,(r,s)\in A^2\mid r+s=0\,\}$.
If we let $Q=\{\,(r,s)\in\A^2\mid r+s=1\,\}$,
then there is a smooth surjective map
$E\times(\A^2\setminus P)\to E\times Q$ such that
$$
\xymatrix{
E\times(\A^2\setminus P) \ar[d] \ar[dr]^{r\varphi+s\psi} \\
E\times Q \ar[r]^{r\varphi+s\psi} & F
}
$$
commutes, and hence the
restriction of $r\varphi+s\psi$ to $E\times Q$ is smooth.
Now let $\lambda\colon R\to Q$ be the map
$t\mapsto t/(t+1)$.
There are induced maps
$$
\xymatrix@C=50pt{
Z_*(F\times R) \ar[r]^(.45){((r\varphi+s\psi)|_Q)^*} &
Z_*(E\times R\times Q) \ar[r]^{\lambda^*} &
Z_*(E\times R\times R)
}
$$
and now by a direct computation we have
$$\partial(\lambda^*((r\varphi+s\psi)|_Q)^*\alpha+\sigma^*\varphi^*\alpha)=
\psi^*\alpha-\varphi^*\alpha$$
for any $\alpha\in Z_*(F\times R)/\pi^*Z_*F$.
\end{proof}

\begin{cor}
There is a functor $X\mapsto \underline A_*X$ and
a natural transformation of functors
$\underline A^\circ_*\to \underline A_*$.
The natural map $\underline A^\circ_*X\to \underline A_*X$
is an isomorphism when $X$ is a scheme.
\end{cor}

\begin{proof}
Proposition \ref{twovbua} implies that for a stack $Y$,
$\{\underline A^\circ_{j+\rk E}E\}$ forms a
direct system over the directed set ${\mathfrak B}_Y$ of
Definition \ref{setby}.
There is, therefore, a functor
$Y\mapsto \underline{\widehat A}_kY:=\varinjlim_{{\mathfrak B}_Y}
\underline A^\circ_{j+\rk E} E$
and a natural map $\underline A_j^\circ Y\to \underline{\widehat A}_jY$.
For projective morphism $f\colon X\to Y$,
we define the groups $\underline {\widehat B}_jX$ as in Definition \ref{bhat},
and then we set
$$\underline A_jY = \varinjlim_{{\mathfrak A}_Y}
( \underline{\widehat A}_jX / \underline{\widehat B}_jX ).$$
As in Remark \ref{almostequiv},
we see that the natural map $\underline A^\circ_*Y\to \underline A_*Y$
is an isomorphism for any scheme $Y$.
\end{proof}

\subsection{The connecting homomorphism}

Let $X$ be a stack, and
let $Y$ be a closed substack of $X$ with complement $U$.  From
the exact sequence of complexes
$$
\xymatrix@C=15pt{
0 \ar[r] &
Z_*(Y\times T^2) \ar[r] \ar[d] &
Z_*(X\times T^2) \ar[r] \ar[d] &
Z_*(U\times T^2) \ar[r] \ar[d] & 0 \\
0 \ar[r] &
Z_*(Y\times R)/\pi^*Z_*Y \ar[r] \ar[d] &
Z_*(X\times R)/\pi^*Z_*X \ar[r] \ar[d] &
Z_*(U\times R)/\pi^*Z_*U \ar[r] \ar[d] & 0 \\
0 \ar[r] &
Z_*Y \ar[r] \ar[d] &
Z_*X \ar[r] \ar[d] &
Z_*U \ar[r] \ar[d] & 0 \\
& 0 & 0 & 0
}
$$
there is a long exact sequence of cohomology groups
\begin{equation}
\label{naiveconnecting}
\underline A^\circ_*Y\to \underline A^\circ_*X \to
\underline A^\circ_*U \stackrel{\delta}\to
A^\circ_*Y \to A^\circ_*X \to A^\circ_*U\to 0.
\end{equation}
A bit of checking shows that the connecting
homomorphism $\delta$ induces a map $\delta$ in the complex
\begin{equation}
\label{connecting}
\underline A_*U\stackrel{\delta}\to A_*Y \to A_*X \to A_*U \to 0.
\end{equation}
The map $\delta$ commutes with flat pullback and projective pushforward.

\begin{pr}
\label{extendedexcision}
Assume $U$ is a global quotient stack.
Then the complex (\ref{connecting}) is exact.
\end{pr}

\begin{proof}
Suppose $h\colon V\to Y$ is projective,
$\alpha\in \hat A_*V$,
and $i_*(h,\alpha)=0$ in $A_*X$,
where $i\colon Y\to X$ denotes inclusion.
After adding components to $V$,
there must exist a projective morphism
$f\colon T\to X$ and projective morphisms $p_1,p_2\colon S\to T$
such that $f\smallcirc p_1$ is 2-isomorphic to $f\smallcirc p_2$,
such that $j\colon V\to Y$ identifies $V$ with $Y\times_XT$,
and such that there exist
$\beta_1\in \hat A^{p_1}_*S$ and
$\beta_2\in \hat A^{p_2}_*S$ such that
$\iota_{p_1}(\beta_1)=\iota_{p_2}(\beta_2)$ in $\hat A_*S$
and
\begin{equation}
\label{alphavanish}
p_{2{*}}\beta_2-p_{1{*}}\beta_1=j_*\alpha.
\end{equation}

Suppose $U\simeq[W/G]$.
Fix a representation bundle $B$ on $BG$.
There exists a projective modification $\pi\colon \widetilde X\to X$
such that some vector bundle $E$ on $\widetilde X$ restricts
to the pullback of $B$ to $U$.
Suppose we are in the special case where $f$
factors as $\pi\smallcirc \tilde f$,
and $\tilde f\smallcirc p_1$ is 2-isomorphic to $\tilde f\smallcirc p_2$.
Then there exists $\beta\in A^\circ_*(\tilde f^*E)^{\oplus n}$
for suitable $n$ such that,
with the fiber diagram
$$
\xymatrix{
R \ar@<-2pt>[d]_{q_1}\ar@<2pt>[d]^{q_2}\ar[r]^k &
S  \ar@<-2pt>[d]_{p_1}\ar@<2pt>[d]^{p_2} \\
V \ar[d]^h \ar[r]^j & T \ar[d]^f \\
Y \ar[r]^i & X
}
$$
there are $\gamma_i\in \hat A^{j\smallcirc q_i}_*R$ such that
$\beta_i=\beta+k_*\gamma_i$ in $\hat A^{p_i}_*S$ for $i=1,2$.
Then, in $A_*X$, we have
$$(h,\alpha)=(h,\alpha+q_{1{*}}\gamma_1-q_{2{*}}\gamma_2)+
(h\smallcirc q_1,\gamma_2-\gamma_1),$$
and each term on the right individually lies in the image of $\delta$.

We now deduce the general case from the special case.
Suppose (\ref{alphavanish}) holds, with $f\colon T\to X$ a general
projective morphism.
There exists
a projective modification $\sigma\colon \widetilde T\to T$ such
that $f\smallcirc \sigma$ factors through $\widetilde X$.
We can make a modification $\tau\colon \widetilde S\to S$ so that
$$(\widetilde S\rightrightarrows \widetilde T\rightarrow \widetilde X)\to
(S\rightrightarrows T\rightarrow X)$$
is a morphism of coequalizer diagrams.
There exists $\tilde \beta_i\in \widehat A^{p_i\smallcirc \tau}_*\widetilde S$
such that
$\tau_*\tilde\beta_i=\beta_i+k_*\gamma_i$
holds in $\widehat A^{p_i}_*S$,
for some $\gamma_i\in \widehat A^{p_i\smallcirc k}R$.
If $\tilde k$ denotes the pullback of $k$,
then we have
$\iota(\tilde\beta_2)-\iota(\tilde\beta_1)=\tilde k_*\varepsilon$
for some $\varepsilon\in \widehat A^{\tau\smallcirc\tilde k}_*\widetilde R$.

Let $\tilde\tau\colon \widetilde R\to R$ denote the
pullback of $\tau$.
By the special case above, we have, in $A_*Y/\delta(\underline A_*U)$,
\begin{align*}
0 &=
(h,\alpha + q_{2{*}}\gamma_2 - q_{1{*}}\gamma_1) -
(h\smallcirc q_1 \smallcirc \tilde\tau, \varepsilon) \\
&=
(h,\alpha) + (h\smallcirc q_1, \gamma_2-\gamma_1-\tilde\tau_*\varepsilon) \\
&= (h,\alpha). \qed
\end{align*}
\renewcommand{\qed}{}\end{proof}

\subsection{Homotopy invariance for vector bundle stacks}

\begin{pr}
\label{fhcgpullback}
Let $X$ be a quasiprojective scheme, let $G$ be a smooth connected
linear algebraic group acting linearly on $X$,
let $Y=[X/G]$, and let $\pi\colon B\to Y$ be a vector bundle stack.
Then $\pi^*\colon \underline A_*Y\to \underline A_*B$ is an isomorphism.
\end{pr}

\begin{proof}
Routine, using scheme approximations (Proposition \ref{schemeapprox} (i)),
plus Proposition \ref{nqiso} and homotopy invariance of $A_*({-};1)$.
\end{proof}

With the machinery developed so far, we finally arrive at a proof
of the homotopy property for vector bundle stacks on
stacks which admit stratifications by global quotients.

\begin{pr}
\label{vbstack}
Let $X$ be a stack which admits a stratification by global quotients,
and let $\pi\colon E\to X$ be a vector bundle stack.
Then $\pi^*\colon A_*X\to A_*E$ is an isomorphism.
\end{pr}

\begin{proof}
We proceed by noetherian induction.
Let $U$ be a nonempty open substack of $X$ such that
$U$ is isomorphic to the quotient of a quasiprojective scheme
by the linearized action of a smooth connected linear algebraic group.

Let $B$ denote the restriction of $E$ to $U$.
By Proposition \ref{schemeapprox} (i), $B$ is a global quotient.
Let $Y=X\setminus U$, and let $F$ denote the restriction of $E$ to $Y$.
Pullback induces a morphism of complexes
$$
\xymatrix{
{\underline A_*B} \ar[r] & A_*F \ar[r] & A_*E \ar[r] & A_*B \ar[r] & 0 \\
{\underline A_*U} \ar[r] \ar[u] & A_*Y \ar[r] \ar[u] &
A_*X \ar[r] \ar[u]^{\pi^*} & A_*U \ar[r] \ar[u] & 0 \ar@{=}[u]
}$$
and by Proposition \ref{extendedexcision} both the bottom complex
and the top complex are exact.  From
Propositions \ref{schemeapprox} (ii) and \ref{fhcgpullback} and
the induction hypothesis,
the five lemma implies that $\pi^*\colon A_*X\to A_*E$ is an isomorphism.
\end{proof}

\section{Intersection theory}

\subsection{Intersections on Artin stacks}
\label{intartsta}

The construction of the deformation space
to a regular local immersion
(section \ref{fmp}) can be generalized to an arbitrary morphism that is
representable and locally separated.
Suppose $f\colon F\to G$ is representable and locally separated.
Then there is a deformation space
$M^\circ_FG\to \bP^1$ with general fiber $G$ and special fiber
the {\em normal cone stack} to $f$.
We will apply this construction to
the diagonal of an Artin stack
to obtain $M^\circ_F (F\times F)\to\bP^1$.
When $F$ is smooth, the special fiber is the tangent bundle stack $TF$,
and this construction leads to an intersection product on $F$ provided that
$F$ admits a stratification by global quotient stacks.

Here is the construction.
Let $f\colon F\to G$ be a representable, locally separated morphism.
It is not hard to see that there exists a commutative square
$$\xymatrix{
U \ar[r] \ar[d] & V \ar[d] \\
F \ar[r] & G
}$$
such that $U$ and $V$ are schemes and such that
the top arrow is a closed immersion.

By the hypotheses on $f$,
the induced map $R:=U\times_FU\to U\times_GV$ is a locally closed immersion,
and is in fact a regular immersion.
We let $S=V\times_GV$, with projections $q_1$ and $q_2$ to $V$.
The composite $R\to U\times_GV\to S$ is also a locally closed immersion,
so there is a deformation space
$M^\circ_RS$, flat over $\bP^1$, with a map
$M^\circ_RS\to S$ that is representable and separated.

By the universal property of blowing up,
the morphism
$M^\circ_RS\to V\times\bP^1$ induced by $q_i$
factors through $\Bl_UV\times\bP^1$, and in fact, factors through
$M^\circ_UV=\Bl_UV\times\bP^1\setminus \Bl_UV$.
So, we get morphisms $r_i\colon M^\circ_RS\to M^\circ_UV$ over $\bP^1$.
for $i=1,2$.
The restriction of $r_i$ over $\bP^1\setminus\{0\}$ is clearly smooth,
and $r_i\times_{\bP^1}\{0\}$ factors as
$$C_RS\to (C_{U\times_VS}S)\times_{U\times_VS}R\simeq C_UV\times_UR\to C_UV$$
and hence is also smooth (the first map appears in the exact sequence
of cones $0\to N_R(U\times_VS)\to C_RS\to C_{U\times_VS}S|_R\to 0$).
The maps $r_i$ are flat by the local criterion for flatness,
and hence are smooth.
The maps $r_i$ together with additional maps obtained in the expected
fashion determine a groupoid
$[M^\circ_RS\rightrightarrows M^\circ_UV]$
which specifies a stack that we denote $M^\circ_FG$.

Now, suppose $f\colon F\to G$ is a representable, locally separated morphism
such that the normal cone stack to $f$
is a vector bundle stack of constant (virtual) rank $d$,
and suppose $G'\to G$ is a morphism
such that $F':=F\times_GG'$ admits
a stratification by global quotient stacks
(this is the case, for instance, if $G'$ admits a stratification by
global quotient stacks).
Then, by the construction of section \ref{fmp} combined with
the homotopy property for vector bundle stacks over a stack which
admits a stratification by global quotient stacks,
we obtain the
{\em refined Gysin homomorphism}
$A_*G'\to A_{*-d}F'$.

The proofs of basic properties about Gysin homomorphisms
(compatibility with flat pullback and projective pushforward,
commutativity, and functoriality) apply unchanged to let us deduce these
properties for the refined Gysin homorphisms appearing in this section.

Let $F$ be a stack which admits a stratification by global quotient stacks.
We remark on several instances when the Gysin homomorphism to
a representable locally separated morphism
$f\colon F\to G$ agrees with the map $f^*$ constructed by
other methods.
First, when $f$ is smooth, or more generally flat and l.c.i., the
Gysin homomorphism is the same as flat pullback:
the virtual normal bundle is a vector bundle stack with
surjective zero section $\varphi\colon {\mathfrak N}\to F$
and the Fulton-MacPherson
produces, starting with the cycle $[Z]$ for some $Z\subset G$,
the cycle $[\varphi^{-1}(f^{-1}(Z))]$ in ${\mathfrak N}$.
Next, when $f$ is a regular local immersion, the virtual normal bundle
is just the usual normal bundle, so the construction reduces to the
usual Fulton-MacPherson construction.
Finally, let $f$ be a l.c.i.\ morphism which admits a global
factorization as a regular immersion $i\colon F\to P$ followed by
a representable locally separated smooth morphism $g\colon P\to G$.
Then there is a morphism of vector bundles
$i^*T_{P/G}\to N_FP$,
and the normal bundle stack to $f$ admits the global presentation
$[N_FP\,/\,i^*T_{P/G}]$.
Now functoriality of the Gysin homomorphism
gives us a new way to see that the definition of
Gysin map to a l.c.i.\ morphism given in \cite{f}, $f^*=i^*\smallcirc g^*$,
is independent of the choice of factorization.
For an l.c.i.\ morphism of schemes which does not admit a global
factorization,
the Gysin map has previously been constructed
in \cite{glms} using higher $K$-theory.

The diagonal morphism of a smooth Artin stack is representable, separated,
and l.c.i.
When $F$ is a smooth Artin stack which admits a stratification
by global quotient stacks, the Gysin homomorphism to the diagonal
$F\to F\times F$ induces a ring structure on $A_*F$.
This is the intersection product.

We can provide an answer to Conjecture 6.6 of \cite{v}.

\begin{thm}
\label{alex}
Let $F$ be a smooth Artin stack which has quasi-finite diagonal,
and let $M$ be a coarse moduli space for $F$.
Then $M$ satisfies Alexander duality.
\end{thm}

\begin{proof}
This is an immediate consequence of the refined Gysin homomorphism
to the diagonal $F\to F\times F$.
In positive characteristic, de Jong's modifications play the role of
resolution of singularity in showing that $\Q$-valued intersection operations
commute.
\end{proof}

\begin{rem}
The notion of Alexander duality in intersection theory was
introduced in \cite{vcomp}.
For a scheme $X$ to satisfy Alexander duality means that $X$ behaves
like a smooth scheme, as far as intersection theory with
rational coefficients is concerned.
In \cite{v} the characteristic zero case of Theorem \ref{alex} is
deduced from the intersection theory on Deligne-Mumford stacks.
\end{rem}

\subsection{Virtual fundamental class}
\label{virfundclass}

Let $X$ be a Deligne-Mumford stack.
The intrinsic normal cone ${\mathfrak C}_X$ is introduced in \cite{bf}
as a tool for constructing a virtual fundamental class in
$A^\circ_*X\otimes\Q$ from a perfect obstruction theory.
A perfect obstruction theory
is an element of the derived category
$E^\smallbullet\in D(\cO_{X_{\rm \acute et}})$
of perfect amplitude contained in $[-1,0]$,
together with a morphism $\varphi\colon E^\smallbullet\to L^\smallbullet_X$,
where $L^\smallbullet_X$ denotes the cotangent complex on $X$,
such that $h^0(\varphi)$ is an isomorphism and
$h^{-1}(\varphi)$ is surjective.
There is an associated geometric object $h^1/h^0(E^\vee)$,
which is a vector bundle stack over $X$:
locally, we can write $E^\smallbullet$ as $[E^{-1}\to E^0]$,
and the quotients $[{E^{-1}}^\vee/{E^0}^\vee]$ patch to form
$h^1/h^0(E^\vee)$.
The map $\varphi$ induces a closed immersion
${\mathfrak C}_X\to h^1/h^0(E^\vee)$.

The construction of the virtual fundamental class proceeds by
starting with the cycle
$[{\mathfrak C}_X]$ on $h^1/h^0(E^\vee)$ and
``intersecting with the zero section'' of
$\pi\colon h^1/h^0(E^\vee)\to X$ to obtain a cycle class on $X$.
As remarked in \cite{bf}, to do this without intersection theory
on stacks requires imposing the additional hypothesis that
$E^\smallbullet$ admits a global presentation as
$[E^{-1}\to E^0]$.
Then we can pull back $[{\mathfrak C}_X]$
to the total space of $E_1:={E^{-1}}^\vee$
and intersect with the zero section
of $E_1\to X$.

The intersection theory of this preprint lets us
remove this extra hypothesis.

\begin{thm}
If $X$ is a Deligne-Mumford stack and
$E^\smallbullet$ is a perfect obstruction theory on $X$
then there exists a unique element $\alpha\in A^\circ_*X\otimes\Q$
such that the pullback of $\alpha$ to $A_*(h^1/h^0(E^\vee))\otimes\Q$
is equal to $[{\mathfrak C}_X]$.
\end{thm}

\begin{proof}
This follows from Propositions \ref{gq} (iii) and \ref{vbstack}.
\end{proof}

\begin{rem}
Using only elementary techniques
(sections \ref{firstproperties} through \ref{elementary}),
it is shown in \cite{krthesis} that
when $X$ is a separated Deligne-Mumford stack,
the map $\pi^*\colon A_*X\otimes\Q\to A_*U\otimes\Q$ is an isomorphism.
The proof uses the fact
(\cite[(10.1)]{l}) that there exists a finite surjective map
from a scheme $T$ to $X$, so by applying Chow's lemma to $T$
we obtain a projective, generically finite, surjective map $f\colon Y\to X$
such that $Y$ is a quasiprojective scheme.
The theorem is true for $Y$
(Proposition \ref{schemeapprox}, with $G$ the trivial group)
and for $Y\times_XY$,
so we deduce the result for $X$
by the analogue for $A_*\otimes\Q$ of the co-sheaf sequences
in $A^\circ_*\otimes\Q$ of \cite{ki}.
\end{rem}

\subsection{Localization formula}
\label{localiz}

We describe a localization formula which is sufficient
for computations in equivariant Chow groups over an algebraically closed
base field of arbitrary characteristic.
We follow \cite{egamj}, which gives a localization formula for
torus actions on schemes.
Using the various functors introduced in this preprint,
we can deduce a similar formula for torus actions on
Deligne-Mumford stacks.
{\em In this section,
all Chow groups are taken to have rational coefficients.}

\begin{lm}
Let $X$ be a stack with quasi-finite stabilizer.
Then $A^\circ_*X\to A_*X$ is an isomorphism and
$\underline A^\circ_*X\to \underline A_*X$ is surjective.
\end{lm}

\begin{proof}
It suffices to show that for any vector bundle $\pi\colon E\to X$
the pullback $A^\circ_*X\to A^\circ_*E$ is an isomorphism and
the map $\underline A^\circ_*X\to \underline A^\circ_*E$ is surjective.
Suppose first that $X$ admits a finite flat cover by a scheme,
$f\colon U\to X$.
Then $f_*\smallcirc f^*$ is an isomorphism,
so we deduce that $\pi$ induces an isomorphism on $A^\circ$ and on
$\underline A^\circ$ from the fact that this holds after pullback via $f$.
By localization, noetherian induction, and the five lemma,
we deduce the desired statements.
\end{proof}

\begin{cor}
\label{fsvanish}
Suppose $X$ has quasi-finite stabilizer.
Then $A_jX=0$ for $j<0$ and
$\underline A_jX=0$ for $j<-1$.
\end{cor}

The localization property for a torus action on a Deligne-Mumford stack
will be a consequence of exactness of
the localization sequence (\ref{connecting}), plus vanishing for
dimension reasons.
When the Deligne-Mumford stack is smooth and has finite stabilizer,
we are able to prove that the fixed stack for the torus action is smooth
with the aid of the next lemma,
and the localization formula follows (the finite stabilizer hypothesis
guarantees that the fixed locus is a closed substack).

\begin{lm}
\label{fixedregular}
Let $A$ be a regular local $k$-algebra with residue field $A/\mm\simeq k$,
let $\sigma$ be a $k$-algebra homomorphism $A\to A$,
and suppose
$\sigma^n=1_A$ for some positive integer $n$,
where $n$ is prime to the characteristic of $k$ in case $k$ has
positive characteristic.
If $I_\sigma$ denotes the ideal generated by $f-\sigma f$ for all $f\in A$,
then $I_\sigma$ is generated by part of a regular sequence of $A$.
\end{lm}

\begin{proof}
Let $d=\dim A$.
The result follows from two facts:
(i) if $f_1$, $\ldots$, $f_d$ is a regular sequence, then
$I_\sigma$ is generated by $(f_1-\sigma f_1)$, $\ldots$,
$(f_d-\sigma f_d)$;
(ii) $f-\frac{1}{n}(f+\sigma f+\cdots+\sigma^{n-1}f) \in I_\sigma$
for any $f\in A$ (we have
$f-\frac{1}{n}(f+\sigma f+\cdots+\sigma^{n-1}f) =
(f-\sigma f)+\frac{n-1}{n}(\sigma f-\sigma^2f) + \cdots +
\frac{1}{n} (\sigma^{n-1}f-f)$).
Using the idempotent projections
$f\mapsto f-\frac{1}{n}\sum_{i=0}^{n-1} \sigma^i f$ and
$f\mapsto \frac{1}{n}\sum_{i=0}^{n-1} \sigma^i f$
we may find a regular sequence $f_1,\ldots, f_d$ such that
$f_j-\frac{1}{n}\sum_{i=0}^{n-1} \sigma^i f_j = f_j$
for $1\le j\le m$
and $\sigma f_j=f_j$ for $m+1\le j\le d$, for appropriate $m$.
Fact (i) identifies $m$ elements which generate $I_\sigma$,
and by (ii) we have $f_j\in I_\sigma$ for $1\le j\le m$.
Since $f_1$, $\ldots$, $f_m$ are linearly independent in
$\mm/\mm^2$, the elements $f_1$, $\ldots$, $f_m$ generate $I_\sigma$.
\end{proof}

Let $G$ be an algebraic group, and let $X$ be a stack.
A $G$-action on $X$ is given by an action map
$a\colon G\times X\to X$ together with an associativity 2-morphism $\gamma$
and identity 2-morphism $\beta$ satisfying compatibility conditions
(associativity must satisfy a commutative cube, and
identity must be compatible with associativity).
Such an action yields a quotient stack $Y$ with map $Y\to BG$
such that $X$ is identified with $Y\times_{BG}\Spec k$.
Suppose the action map is the projection map $pr_2\colon G\times X\to X$.
Then the associativity 2-morphism is an automorphism
$pr_3\colon G\times G\times X\to X$,
and the action is trivial (i.e., $Y\simeq X\times BG$) if and only if
$\gamma$ is obtained from $1_{pr_3}$ by composing with $\beta$.
A general action is trivial if there is a 2-morphism $\delta\colon a\to pr_2$
such that
if $\gamma'$ and $\beta'$ denote the compatibility morphisms obtained
from $\gamma$ and $\beta$ respectively by applying $\delta$,
then $\gamma'$ is obtained from $1_{pr_3}$ by composing with $\beta'$.
So, for example, if $G$ is connected, then the $G$-action is trivial
if and only if $a$ is 2-isomorphic to $pr_2$.
If we form the stabilizer diagram
\begin{equation}
\label{stdiag}
\begin{split}
\xymatrix{
S\ar[r] \ar[d]_\varphi & X \ar[d]^{diag_X} \\
G\times X \ar[r]^{(pr_2,a)} & X\times X
}
\end{split}
\end{equation}
then $S$ has the structure of group object in the category of
$X$-stacks, $\varphi$ is
a group homomorphism, and if $G$ is connected then the action is trivial
if and only if $\varphi$ admits a splitting.

{\em For the remainder of this section, we assume the base field $k$ to
be algebraically closed.}

Suppose now $T$ is the one-dimensional torus $T={\mathbb G}_m$, and
suppose $T$ acts on $X=BH$, for some finite group $H$.
We form the stabilizer diagram;
the coarse moduli space to the component of the
identity of the stabilizer $S$ is some $T'\simeq {\mathbb G}_m$,
with $T'\to T$ unramified.
If we consider the induced $T'$-action on $X$, then
$\varphi$ in (\ref{stdiag}) admits a splitting, and thus the action is trivial.
Requiring a torus extension to
trivialize the action on fixpoints accounts for the fractional weights
that show up in calculations, e.g., in \cite{ko}.

\begin{pr}
\label{fixedpoint}
Let $X$ be an integral
Deligne-Mumford stack, and assume $X$ has finite stabilizer.
Suppose the one-dimensional torus $T$ acts on $X$.
Then there exists a $T$-stable closed substack $Z$ of $X$ and a
positive integer
$n$, prime to the characteristic of the base field,
such that if $T'\to T$ denotes the $n$-fold cover,
then the induced action of $T'$ on $Z$ is trivial,
and the induced action of $T'$ on $X\setminus Z$ has quasi-finite
stabilizer.
Moreover, if $X$ is smooth, then $Z$ is smooth.
\end{pr}

\begin{proof}
Consider the stabilizer diagram (\ref{stdiag}), with $G=T$.
Since $X$ has finite stabilizer, the fiber of $\varphi$ over
a general point of $T$ maps finitely to $X$, and determines a
closed substack $Z$ of $X$ (which we give the reduced substack structure).
For suitable $n$, the $n$-fold covering torus $T'$ acts trivially
on the integral zero-dimensional substacks of $Z$.
Consider now diagram (\ref{stdiag}) with $G=T'$ and $X=Z$.
The map from the connected component of the identity to $T'\times Z$
is an isomorphism, and hence $T'$ acts trivially on $Z$.
The action of $T'$ on $X\setminus Z$
clearly has quasi-finite stabilizer.

Now assume $X$ is smooth.
We find $Z=\varphi(\varphi^{-1}(\{t\}\times X))$ for all $t$ outside
a finite subset of $T$, and in particular for $t$ equal to
a primitive $r^{\rm th}$ root of unity, for suitable $r$
prime to the characteristic of the base field.
So, $Z$ is the fixed locus of a cyclic group action.
The cyclic group action may be presented
by the action of a generator $\sigma\colon X\to X$,
together with a 2-morphism
$\delta\colon \sigma^r\to 1_X$.
Since the action on $Z$ is trivial, there exists
a 2-morphism $\sigma|_Z\to 1_Z$,
compatible with $\delta$.

Let $f\colon U\to X$ be an \'etale atlas,
and let $R=U\times_XU$.
Replacing $U$ by the fiber product of $f$, $f\smallcirc \sigma$, $\ldots$,
we may assume that $\sigma$ is represented by an automorphism of $U$.
The closed substack $Z$ of $X$ corresponds to a closed subscheme $Y$ of $U$.
Suppose $u$ is a closed point of $Y$.
Passing to the henselization $U^h$ of $U$ at $u$,
we may assume $\sigma|_{Y^h}=1_{Y^h}$.
Thus we have $\sigma^r=1_{U^h}$.
Hence the ideal $I$ of $\cO_{U^h}$
corresponding to the closed subscheme $Y$ of $U$
satisfies $I\supset I_\sigma$, where $I_\sigma$
corresponds to the fixed locus of the
action of $\sigma$ on $U^h$.
Since $R$ contains an identity component, we also have
$I\subset I_\sigma$.
So, by Lemma \ref{fixedregular}, $\cO_{u,Y}$ is regular, and thus
$Z$ is smooth.
\end{proof}

\begin{thm}
Let $X$ be a Deligne-Mumford stack which has finite stabilizer,
and suppose $T={\mathbb G}_m$
acts on $X$.
Let $A^T_*X$ denote $A_*$ of the stack quotient $[X/T]$.
If $X^T$ denotes the fixed-point substack (which was called $Z$ in
Proposition \ref{fixedpoint}), and if we let
$t=c_1(\cO(1))\in A_*BT$, then
the inclusion $Z\to X$ induces an isomorphism
$$A_*X^T\otimes_{\Q} \Q[t,t^{-1}]\to A^T_*X\otimes_{\Q[t]} \Q[t,t^{-1}].$$
\end{thm}

\begin{proof}
Let $T'\to T$ be the $n$-fold covering torus from Proposition \ref{fixedpoint}.
The action of $T'$ on $U=X\setminus X^T$ has quasi-finite stabilizer,
so the zeroth and first Chow homology groups of the stack quotient
$A^T_*U$ and $\underline A^T_*U$ vanish in the dimension ranges
indicated by Corollary \ref{fsvanish}.

It suffices to prove the theorem for the induced $T'$-action on $X$,
so we may assume $T$ acts trivially on $X^T$.
Let $i\colon X^T\to X$ denote inclusion.
Exactness of the localization sequence
$$\underline A^T_jU \to A_jX^T\otimes_{\Q}\Q[t]\stackrel{i_*}\to
A^T_jX\to A^T_jU\to 0$$
implies that both the kernel and cokernel of $i_*$ are killed by
suitable powers of $t$, and the theorem follows.
\end{proof}

Let $X_1$, $\ldots$, $X_m$ denote the connected components of $X^T$,
and let $j_i$ denote the corresponding inclusion maps, $j=1$, $\ldots$, $m$.
Apparently, any $\alpha\in A^T_*X\otimes_{\Q[t]}\Q[t,t^{-1}]$
can be written uniquely as
$\alpha=\sum_j i_{j{*}}\alpha_j$
with $\alpha_j\in A_*X_j\otimes\Q[t,t^{-1}]$.
Now suppose $X$ is smooth.
In analogy with the construction of \cite{eg} we obtain a ring structure
on $A^T_*X$.
By Proposition \ref{fixedpoint}, each $X_i$ is smooth,
so $A^T_*X_i\simeq A_*X_i\otimes\Q[t,t^{-1}]$ also has a ring structure.

Writing $1=\sum_j i_{j{*}}\omega_j$ in $A^T_*X$,
the identity $\alpha=\alpha\cdot 1$
implies the identities
$\alpha_j=\omega_j\cdot i^*_j\alpha$ in the ring $A_*X_j\otimes\Q[t,t^{-1}]$,
for each $j$.
In particular,
$1=\omega_j\cdot c_{\rm top}(N_{X_j}X)$,
so we deduce

\begin{cor}
Assume $X$ is a smooth Deligne-Mumford stack which has
finite stabilizer,
and suppose the one-dimensional torus $T$ acts on $X$
with fixed locus $X^T=X_1\amalg\cdots\amalg X_m$.
Then the equivariant normal bundle $c_{\rm top}(N_{X_j}X)$ is invertible in
$A_*X_j\otimes\Q[t,t^{-1}]$ for each $j$, and we have
$$\alpha=\sum_{j=1}^m i_{j{*}}
\frac{ i^*_j\alpha }{ c_{\rm top}(N_{X_j}X) }$$
for any $\alpha\in A^T_*X\otimes_{\Q[t]}\Q[t,t^{-1}]$.
\end{cor}

\noindent
Department of Mathematics\\
University of Pennsylvania\\
209 South 33rd Street\\
Philadelphia, PA 19104
\end{document}